\documentclass[12pt,reqno]{amsart}
\usepackage{amsthm,amsmath,amssymb,amsfonts,amscd,latexsym,verbatim,xypic}
\def\suffix{ps}
\newcount\system
\global\system=3   % for unix(dvips)
\def\ifundefined#1{\expandafter\ifx\csname#1\endcsname\relax}
\ifundefined{figdir}\def\figdir{}\fi
\newcount\firstline
\newdimen\pswidth  \newdimen\xleft
\newdimen\psheight \newdimen\ytop \newdimen\ybot
\newcount\justx \newcount\justy
\global\justx=0 \global\justy=0
\newdimen\vpos \newtoks\labeL 
\newread\labeLfile \newdimen\xcoord \newdimen\ycoord
\newif\ifdoit 
\newbox\labox
\newdimen\xdvikwid 
\newdimen\xdvikht
\newdimen\pspoints
\newdimen\rwi
\newcount\temp
\def\readdim#1{\global\read\labeLfile to \temp
\global #1=\temp pt}
\def\figcrop#1{\par%  #1=filename
\openin\labeLfile=\figdir#1.lbl                                              
\global\read\labeLfile to\firstline\message{#1}               
\global\read\labeLfile to\temp%read overall dimensions                                     
\readdim{\ybot}
\readdim{\xleft}%               read upper left point
\readdim{\ytop}
\global\read\labeLfile to\justx%ignore
\global\read\labeLfile to\justy%ignore
\global\read\labeLfile to\labeL%ignore
\readdim{\pswidth}%            read lower right point
\global\advance\pswidth by -\xleft
\readdim{\psheight}
\global\advance\ybot by -\psheight
\global\advance\psheight by -\ytop
\global\read\labeLfile to\justx%ignore
\global\read\labeLfile to\justy%ignore
\global\read\labeLfile to\labeL%ignore                                    
\vbox to\psheight{\vfill
\ifnum\system=1
\fi %textures
\ifnum\system=2
\fi %msdos
\ifnum\system=3
                                                 \fi         %%unix:dvips
\ifnum\system=4
\fi         %%unix:dvips,scaled
\ifnum\system=1
\hbox to \pswidth{\kern-\xleft\special{postscriptfile \figdir#1.\suffix }\hfil}\fi
\ifnum\system=2
\hbox to \pswidth{\kern-\xleft\special{ps: plotfile \figdir#1.\suffix }\hfil}\fi
\ifnum\system=3
\hbox to \pswidth{\kern-\xleft\includegraphics{\figdir#1.\suffix}\hfil}\fi
\ifnum\system=4
\hbox to \pswidth{\kern-\xleft\includegraphics{\figdir#1.\suffix}\hfil}\fi
\ifnum\system=5
\hbox to \pswidth{\kern-\xleft\includegraphics{\figdir#1.\suffix}\hfil}\fi %orphee
\ifnum\system=6
   \xdvikwid=\pswidth
   \xdvikht=\psheight
   {\global\divide\xdvikwid by \pspoints}
   {\global\divide\xdvikht by \pspoints}
   \rwi=\xdvikwid
    {\global\multiply\rwi by 10}
\hbox to \pswidth{\kern-\xleft\includegraphics{\figdir#1.\suffix\space}\hfil}\fi                   %xdvik
\vskip -\baselineskip
\vskip -\ybot 
\vskip-\psheight %                                     
\hbox to\pswidth  {\hss}%                                            
\parindent=0pt\offinterlineskip                                       
\vpos=0 pt%                                                              
\loop\readdim{\xcoord}                                 
\ifdim \xcoord < -999pt \doitfalse\else\doittrue\fi                        
\ifdoit \advance \xcoord by -\xleft
\readdim{\ycoord}
\advance \ycoord by -\ytop                              
\global\read\labeLfile to\justx                                       
\global\read\labeLfile to\justy                                       
\global\read\labeLfile to\labeL
\global\setbox\labox=\hbox{\labeL\hskip-0.3em}%    
\advance\vpos by-\ycoord                                              
\vskip-\vpos \vpos=\ycoord                                         
\hbox to\pswidth{\hskip\xcoord %                                 
\hbox to 0pt{\ifnum\justx>0\hss\fi%                                   
\vbox to0pt{%                                                         
\ifnum\justy<2\vss\fi%                                                
\copy\labox\kern0pt%  
\ifnum\justy>0\vss\fi}%                                               
\ifnum\justx<2\hss\fi}%                                               
\hss}%                                                                
\repeat%                                                              
\advance\vpos by-\psheight%                                           
\vskip-\vpos %                                                     
}\closein\labeLfile}
\def\figplace#1#2#3{
\openin\labeLfile=\figdir#1.lbl
\ifeof \labeLfile
       \immediate\write16{***Can't find \figdir#1.lbl; Skipping it.***}
\else  \closein\labeLfile
       \null\hskip#2\raise #3 \hbox{\figcrop{#1}}
\fi
}
\def\bbbone{{\mathchoice {\rm 1\mskip-4mu l} {\rm 1\mskip-4mu l}
{\rm 1\mskip-4.5mu l} {\rm 1\mskip-5mu l}}}

\newcommand{\be} {\begin{equation}}
\newcommand{\ee} {\end{equation}}
\newcommand{\bea} {\begin{eqnarray}}
\newcommand{\eea} {\end{eqnarray}}
\newcommand{\lp}  {\left(}
\newcommand{\rp}  {\right)}
\newcommand{\cU}{\mathcal U}

\newcommand{\cN}{\mathcal N}
\newcommand{\cP}{\mathcal P}

\newcommand{\cS}{\mathcal S}
\newcommand{\cW}{\mathcal W}

\newcommand{\al}{\alpha}

\newtheorem{Theorem}{Theorem}[section]
\newtheorem{Lemma}[Theorem]{Lemma}
\newtheorem{Proposition}[Theorem]{Proposition}
\newtheorem{Corollary}[Theorem]{Corollary}
\newtheorem{Remark}[Theorem]{Remark}

\newcommand{\complex}{\mathbf C}
\newcommand{\Z}{\mathbf Z}

\renewcommand{\P}{\mathbb P}        
\newcommand{\T}{\mathbb T}        
\renewcommand{\O}{\mathcal O} 
\newcommand{\C}{\mathcal C} 
\newcommand{\gC}{\mathfrak C} 
\newcommand{\I}{\mathcal I} 
\newcommand{\J}{\mathcal J} 
\newcommand{\Q}{\mathcal Q} 
\newcommand{\M}{\mathcal M} 
 
\newcommand{\R}{\mathfrak P}

\newcommand{\Low}{{\mathcal L}\,}
\newcommand{\Hi}{{\mathcal H}\,}
\newcommand{\E}{\mathcal E} 
 
\newcommand{\ux}{\mathbf x} 
\newcommand{\uy}{\mathbf y} 
\newcommand{\uz}{\mathbf z} 
\newcommand{\A}{\mathcal A}  
\newcommand{\U}{\mathcal U}  
\newcommand{\D}{\mathfrak D}
\newcommand{\He}{\text{He}}
\newcommand{\period}{\, .}
\newcommand{\Gordan}[9]{ 
\left(\begin{array}{ccc} #1 & #2 & #3 \\ #4 & #5 & #6 \\ #7 & #8 & #9
\end{array} \right)}
\newcommand{\ra}{\rightarrow}

\newcommand{\lra}{\longrightarrow}

\renewcommand{\ker}{\text{ker}\,}
\newcommand{\coker}{\text{coker}\,}

\newcommand{\spec}{\text{Spec}\,}
\newcommand{\demo}{\noindent {\sc Proof.}\;}
 
\begin{document}
\title{The Bipartite Brill--Gordan Locus and angular momentum} 
\author[Abdesselam and Chipalkatti]
{Abdelmalek Abdesselam and Jaydeep Chipalkatti} 
\maketitle 

\parbox{12cm}{\small 
{\sc Abstract.} 
Given integers $n,d,e$ with $1 \le e < \frac{d}{2}$, let 
$X \subseteq \P^{\binom{d+n}{d}-1}$ denote the locus of 
degree $d$ hypersurfaces in $\P^n$ which are supported on two 
hyperplanes with multiplicities $d-e$ and $e$. Thus $X$ is the 
{\sl Brill-Gordan} locus associated to the partition $(d-e,e)$. 
The main result of the paper is an exact determination 
of the Castelnuovo regularity of the ideal of $X$. Moreover 
we show that $X$ is $r$-normal for $r \ge 3$.}

\medskip 

\parbox{12cm}{\small 
In the case of binary forms (i.e., for $n=1$) we give an invariant theoretic 
description of the ideal generators, and furthermore exhibit 
a set of two covariants which define this locus 
set-theoretically.}

\medskip 

\parbox{12cm}{\small 
In addition to the standard cohomological tools in algebraic 
geometry, the proof crucially relies on the nonvanishing
of certain $3j$-symbols from the quantum 
theory of angular momentum.}

\vspace{5mm} 

\mbox{\small AMS subject classification (2000): 
14F17, 20G05, 22E70, 33C20.}

\bigskip 

\parbox{12cm} 

%%%%%%%%%%%%%%%%%%%%%%%%%%%%%%%%%%%%%%%%%%%%%%%%%%%%%%%

\bigskip

\section{Introduction} 
This paper is a sequel to \cite{AC1}, to which we refer the reader for a 
detailed introduction to the problem considered here. However it may 
be read by itself without substantial loss of continuity. 

\subsection{} \label{section.notation} 
The base field will be $\complex$. Let $V$ denote an $(n+1)$-dimensional 
complex vector space, with 
$W = V^* = \text{span} \, \{x_0,x_1,\dots,x_n \}$. 
The degree $d$ homogeneous forms in the $x_i$ 
(distinguished up to scalars) are parametrized by 
\[ \P^N = \P^{\binom{d+n}{d}-1} = \P \, S_d \, W = 
\text{Proj} \, R, \] 
where $R$ is the symmetric algebra
$\bigoplus\limits_{r \ge 0} S_r(S_d \, V)$. Now let $e$ be an 
integer such that $1 \le e \le \frac{d}{2}$, and define 
the $2n$-dimensional subvariety 
\[  X^{(d-e,e)} = \{ [F] \in \P^N : F = L_1^{d-e} \, L_2^e \; \text{for some 
linear forms $L_1,L_2$} \}. \] 
We will merely write $X$ for $X^{(d-e,e)}$ if no confusion is likely. 
In the language of \cite[\S 1]{AC1}, this is the {\sl Brill-Gordan locus} 
associated to the partition $(d-e,e)$. In the 1890s, Brill and Gordan 
considered the problem of finding defining equations for the following 
variety 
\[ \{ [F] \in \P^N : F = \prod\limits_{i=1}^d  L_i 
\; \; \text{for some linear forms $L_i$} \}, \] 
this serves as the motivation behind this nomenclature. 

Now assume $ e < \frac{d}{2}$ (the case $e = \frac{d}{2}$ was treated 
in \cite{AC1}), and consider the graded ideal $I_X \subseteq R$. 
The main result of this paper is the following: 

\begin{Theorem}\sl 
The Castelnuovo regularity of the ideal $I_X$ is equal to 
\begin{equation}
m_0 =\lceil \max \, \{4, \, n+2 + \frac{1-n}{d}, \,
2n+1 - \frac{n}{e}\}\rceil. 
\label{defn.m0} \end{equation} 
\emph{A fortiori}, the ideal is generated by forms of degree 
at most $m_0$. 
\label{main.theorem} \end{Theorem} 
During the course of the proof the following result emerges naturally. 
\begin{Proposition} \sl 
For $r \ge 3$, the variety $X$ is $r$-normal, i.e., the morphism 
\[ H^0(\P^N,\O_{\P^N}(r)) \lra H^0(\P^N,\O_X(r))  \] 
is surjective. 
\end{Proposition} 

The imbedding $X \subseteq \P^N$ is stable for the natural action of 
the group $SL(V)$; this fact is essentially used throughout the 
paper. 
\subsection{Binary Forms} 
A particularly interesting case 
is that of binary forms (i.e., $n=1$), when we get  
$m_0 = 4$. Together with \cite[Theorem 1.4]{AC1}, this completely
proves the following result which was first conjectured in 
\cite{Chipalkatti1}. 

\begin{Theorem} \sl With notation as above, 
the Castelnuovo regularity of $I_X$ is equal to $3$ if $2 \, e = d$, 
and $4$ otherwise. 
\end{Theorem} 

Let $m_0$ stand for either $3$ or $4$. 
An irreducible $SL_2$-submodule 
\[ S_q \subseteq (I_X)_{m_0} \subseteq S_{m_0}(S_d), \] 
corresponds to a covariant of binary $d$-ics 
of degree $m_0$ and order $q$ which identically vanishes on $X$. 
For the case $d = 2 \, e$, we had explicitly described all such 
covariants in \cite[\S 7]{AC1} as linear combinations of 
compound transvectants. In Section \ref{binary.idealgen} below
we outline such a description 
for the case $d \neq 2 \, e$. However, in this case the expressions
are not as explicit as before, to the extent that they involve 
{\sl Clebsch-Gordan coefficients}. 
In general, no `closed formulae' are known for the latter. 

Section~\ref{section.binary.set-th} is independent of the rest of the 
paper. There we construct two covariants $\D$ and $\gC_e$ which define 
the locus 
$X^{(d-e,e)}$ set-theoretically, i.e., for a binary $d$-ic $F$, 
\[ F \in X^{(d-e,e)} \iff \gC_e(F) = \D(F) = 0.\] 
The constructions are elementary, in fact they involve little beyond the 
Hessian and the Wronskian determinants. 

\subsection{} 
We begin proving the main theorem in Section~\ref{section.conductorsheaf}. 
Let $\I_X \subseteq \O_{\P^N}$ denote the ideal sheaf of $X$, 
then the statement to be established is 
\[ H^q(\P^N,\I_X(m_0-q)) = 0 \quad \text{for $q \ge 1$.}  \] 
First we determine the `conductor sheaf' supported on the singular locus 
of $X$, and then calculate its cohomology using the Borel-Weil-Bott theorem; 
this gives the required vanishing for $q \ge 2$. 
The case $q=1$ occupies the bulk of the paper. We reduce  it to a problem 
about transvectants of binary forms, and then settle the latter using 
some explicit combinatorial calculations in 
Sections \ref{section.trans} through \ref{section.proof.I.II}. 

As mentioned earlier, Sections~\ref{section.binary.set-th} 
and \ref{binary.idealgen} are devoted to 
binary forms. In Section~\ref{x32} we treat the variety $X^{(3,2)}$ for 
ternary quintics. Using some machine calculations (in Macaulay-2), we 
express the ideal generators of $X$ as concomitants of ternary quintics. 

The basic representation theory of $SL(V)$ may be found in~\cite{FH}. 
All the terminology from algebraic geometry agrees with~\cite{Ha}. 
As for the classical invariant theory and the symbolic method, 
\cite{Glenn, GrYo} will serve as our standard references.
Note however that the {\sl interpretation} of the symbolic method
which we will use is the one briefly given in \cite[\S 1.7]{AC1}.
For the benefit of the reader with no prior familiarity with this
somewhat controversial tool, 
a more detailed explanation would be appropriate here; it is
provided in the following section.

\section{The classical symbolic method}
\label{section.symbolic.method}
The symbolic notation was introduced by
Aronhold in~\cite{Aronhold}, and most prominently 
developed by the German school of invariant theory led by 
Clebsch and Gordan~\cite{Clebsch, Gordan}.
It is in fact a powerful reformulation
of Cayley's theory of hyperdeterminants~\cite{Cayley}.
Most presentations of this method go as follows.

\subsection{} Let
\[ A(x_0,x_1)=\sum_{i=0}^m \, 
\binom{m}{i} \alpha_i\ x_0^{m-i} x_1^i \]
be a generic binary form of degree $m$.
Write $A$ `symbolically' as
\[ A(x_0,x_1)=(a_0 \, x_0+a_1 \, x_1)^m, 
\]
i.e., one `postulates' that 
\begin{equation} 
\alpha_i=a_0^{m-i} a_1^i \quad 
\text{for $0 \le i \le m$.} \label{postul} \end{equation} 
However, this would introduce unwanted relations such as 
$\alpha_0 \, \alpha_2=\alpha_1^2$. In order to prevent these, 
the prescription is to use different symbols for each
individual factor $\alpha_i$ in a product like $\alpha_0 \, \alpha_2$.
One therefore introduces additional letters and writes
\[
A(x_0,x_1)=(a_0 \, x_0+a_1 \, x_1)^m=(b_0 \, x_0+b_1 \, x_1)^m= \dots, 
\]
and then the translation between monomials in the coefficients of $A$
and those in the symbolical letters becomes 
\[
\alpha_0 \, \alpha_2=a_0^m \, b_0^{m-2} \, b_1^2, \quad 
\alpha_1^2 \, \alpha_2=
a_0^{m-1} \, a_1 \, b_0^{m-1} \, b_1 \, c_0^{m-2} \, c_1^2, 
\; \text{etc}. 
\]
Needless to say, this explanation is far from satisfactory
and could understandably seem, on a first encounter, closer
to witchcraft than mathematics. In the recent mathematical 
literature confronting this issue, one can trace essentially
two different attitudes towards the symbolic method.
The first one is simply to ignore it altogether and do without
it completely; however this entails throwing away a very valuable
tool and comes at a cost: missing some of the gems of classical
geometry which are given an appropriate display, for instance 
in~\cite{Hunt}. The second is the compromise 
expressed in [loc.~cit.~pp.~290-291], where the use of this method 
is advocated regardless of rigor as a quick way to {\sl guess} 
polynomial identities involving invariants; while 
the task of {\sl checking} these identities is left for other methods, 
for instance the help of a computer. 

This is somewhat analogous to the situation with the recent 
cross-fertilization between algebraic geometry and 
theoretical physics (see e.g.~\cite{Hori}).
In a few but important instances, mathematical statements 
heuristically derived by physicists using functional integral methods 
were later established as theorems, either by using previously 
existing tools of algebraic geometry, or by devising new ones
in order to bypass path integrals.
In the latter situation, there is indeed the genuine difficulty
of making functional integration rigorous,
which is the business of constructive field theory
(see e.g.~\cite{Gawedzki}).
However, in the case of the classical symbolic method of invariant theory,
with only basic multivariate calculus as a prerequisite, we will
show that there is {\sl no difficulty at all}.

\subsection{} 
The simple trick is to use the easily checked identity
\[ A(\ux)=\frac{1}{m!} \, A(\frac{\partial}{\partial a_0},
\frac{\partial}{\partial a_1})
\ a_\ux^m
\]
where $a_\ux=a_0 \, x_0+ a_1 \, x_1$. By convention, the differential
operators apply to whatever is on the right, and
this equation, as well as the ones that follow, are to be 
understood {\sl verbatim} rather than `interpreted symbolically'.
We will use the notation $\int_A {\rm d}a$ for the differential operator
\[
\frac{1}{m!} \, A(\frac{\partial}{\partial a_0},
\frac{\partial}{\partial a_1}), 
\]
so that the previous equation becomes
\[
A(\ux)=\int_A {\rm d}a\ a_\ux^m\ .
\]
This choice of notation can be justified by the following
reasons :

\noindent
1.~In practice, manipulating symbolic letters in the same
way as dummy variables of integration is enough to guard 
against computational blunders.

\noindent
2.~An alternate way to put the symbolic method on a rigorous
footing, given by Littlewood~\cite[pp.~326--327]{Littlewood},
precisely uses the fact that any form can be written
as the sum of sufficiently many powers of linear forms.
The idea is similar in spirit to the use of one's
favourite integral representation, such as the Fourier transform,
for doing analytic calculations.

\noindent
3.~Last but not least, it takes less room on the page.

\subsection{} 
Now consider two binary forms $A(\ux),B(\ux)$ of
respective degrees $m,n$, and let $k$ be an integer such that 
$0\le k\le \min\{m,n\}$. By definition, 
the $k$-{\sl{th transvectant}} of $A$ and $B$ is the degree 
$m+n-2k$ form given by
\begin{equation} 
(A,B)_k(\ux)= \frac{(m-k)! \, (n-k)!}{m! \, n!} 
\left.
\left\{ \Omega_{\ux\uy}^k \, 
A(\ux) \, B(\uy) \right\}
\right|_{\uy:=\ux} \label{transvectdef} 
\end{equation} 
where
\[
\Omega_{\ux\uy}=\frac{\partial^2}{\partial x_0 \, \partial y_1}-
\frac{\partial^2}{\partial x_1 \, \partial y_0}
\]
is Cayley's Omega operator, and $\uy=(y_0,y_1)$ is an extra set of variables.
Alternately one can also expand $\Omega_{\ux\uy}$ by the binomial
theorem, and write the equally useful formula
\begin{equation} 
(A,B)_k(\ux)=
\frac{(m-k)! \, (n-k)!}{m! \, n!} 
\sum\limits_{i=0}^k (-1)^i \binom{k}{i} \, 
\frac{\partial^k A}{\partial x_0^{k-i} \, \partial x_1^i} 
\frac{\partial^k B}{\partial x_0^i \, \partial x_1^{k-i}}\ .
\label{trans.formula} 
\end{equation} 
Now we have 
\[ \begin{aligned}
(A,B)_k(\ux) & = 
\frac{(m-k)! \, (n-k)!}{m! \, n!}
\left. \{ 
\Omega_{\ux\uy}^k \, (\int_A {\rm d}a\ a_\ux^m) 
(\int_B {\rm d}b\ b_\uy^n) \}
\right|_{\uy:=\ux} \\
& = \frac{(m-k)! \, (n-k)!}{m!\, n!}
\int_A {\rm d}a\ \int_B {\rm d}b\ 
\left[ \left.
\{ \Omega_{\ux\uy}^k \, 
a_\ux^m \, b_\uy^n \}\right|_{\uy:=\ux} \right]; 
\end{aligned} \]
because the differential operators
$\int_A {\rm d}a, \int_B {\rm d}b$ {\sl commute} with
$\Omega_{\ux\uy}$ and substitution of $\ux$ into $\uy$, 
for the mere reason that the variables $a$, $b$, $\ux$, $\uy$
are distinct.
Now an elementary calculation using Leibnitz's rule shows that
\[ \left. 
\{ \Omega_{\ux\uy}^k \, a_\ux^m \, b_\uy^n \} \right|_{\uy:=\ux} =
\frac{m! \, n!}{(m-k)! \, (n-k)!} \, 
(a \, b)^k \, a_\ux^{m-k} \, b_\ux^{n-k} 
\]
where the so-called symbolic bracket $(a \, b)$ is a 
compact notation for $a_0 \, b_1-a_1 \, b_0$. Therefore
\[ (A,B)_k=\int_A {\rm d}a \int_B {\rm d}b \; I(a,b,\ux), 
\]
where the `integrand' $I(a,b,\ux)$, i.e., {\sl the classical symbolical
expression} for the transvectant
$(A,B)_k$, is equal to $(a \, b)^k \, a_\ux^{m-k} \, b_\ux^{n-k}$.
The particularly nice final expression
justifies the combinatorial normalization factor in the original
definition (\ref{transvectdef}). A modern physicist might
say that the classical mathematicians had the consummate wisdom
of normalizing `sums over Wick contractions' as probabilistic
{\sl averages}.
One should bear in mind that as long as the forms $A,B$
are generic, it is preferable to do the calculations
on the `integrand' and refrain from actually performing
the `integral'. However, this no longer applies as soon as one 
substitutes a composite algebraic expression for one of the forms,
for instance a decomposition into linear factors involving the roots
(as in~\cite[p.~18]{AC1}), or the transvectant of two other
forms etc.

Here is not the place for a comprehensive review (as yet unwritten) 
of the classical symbolic method, 
and its relation to the calculus of Feynman diagrams
(see~\cite{AC1}) as well as the quantum theory of angular momentum
(see~\cite{Biedenharn, Biedenharn2}).
Nevertheless, the above should suffice in order to enable
the reader to check the calculations in 
Sections~\ref{section.trans}~through~\ref{section.proof.I.II} 
with all the necessary mathematical rigor.

\section{The conductor sheaf} \label{section.conductorsheaf}
Let us pick up the thread from the beginning of 
section~\ref{section.notation}. We have an  $n$-dimensional 
smooth subvariety 
\begin{equation} 
Z = \{ [F] \in \P^N: F = L^d \;\; \text{for 
some $L \in W$} \} \subseteq X, \end{equation} 
which is the $d$-fold Veronese imbedding of $\P^n = \P W$ into $\P^N$. 
There is a proper birational morphism 
\begin{equation} \P W \times \P W \stackrel{f}{\lra} \P S_d \, W, \quad 
(L_1,L_2) \lra L_1^{d-e} \, L_2^e
\end{equation}
with image $X$. It is an isomorphism 
over $X \setminus Z$, hence we have an exact sequence 
\begin{equation}
0 \lra \O_X \lra f_* \O_{\P W \times \P W} \lra \Q \lra 0,
\label{conductor.seq} \end{equation}
where the {\sl conductor sheaf} $\Q$ 
(so called because $f$ is the normalization of $X$) has support $Z$. 
Let $\delta: \P W \lra \P W \times \P W$ 
denote the diagonal imbedding, and $g = f \circ \delta$. 
\[ \diagram 
\P W \dto_{\delta} \drto^g \\ 
\P W \times \P W \rto^f & \P S_d \, W 
\enddiagram \] 
\begin{Proposition} \sl 
The pullback $g^* \Q$ is isomorphic to $\Omega^1_{\P W}$ (the cotangent sheaf of $\P W$). 
\end{Proposition}
\demo 
Let $\J$ denote the ideal sheaf of image$(\delta)$, then we have an 
exact sequence 
\[ 0 \lra \J \lra \O_{\P W \times \P W} \lra 
\delta_* \O_{\P W} \lra 0. \] 
{\bf Claim 1}: The inclusion $\J \subseteq \O_{\P^n \times \P^n}$ 
factors through the natural map (see \cite[p.~110]{Ha})
\[ f^*f_* \O_{\P^n \times \P^n} \lra \O_{\P^n \times \P^n}. \]
Proof. 
Since $f$ is an affine morphism, the claim is local on $\P^n \times \P^n$. 
Hence, restricting to affine open sets, we may write 
\[ f: \spec B \lra \spec A, \quad 
  \delta: \spec B/J \lra \spec B. \] 
Then $f^*f_* \O_{\P^n \times \P^n}$ is locally represented 
by the $B$-module $B \otimes_A B_A$ (where $B_A$ denotes $B$ considered 
as an $A$-module). The $B$-module map 
\[ J \lra B \otimes_A B_A, \quad x \lra x \otimes 1 \] 
is the required factorization, which proves Claim 1. 
Composing with 
\[ f^*f_* \O_{\P^n \times \P^n} \lra f^* \Q, \] 
we get a map $\J \stackrel{q}{\lra} f^*\Q$. 

\medskip 

\noindent {\bf Claim 2}: $q$ is surjective.
\smallskip 

\noindent Proof. It will suffice to show that the composite 
$f^*f_* \J \lra \J \stackrel{q}{\lra} f^*\Q$ is surjective. 
We have a commutative ladder 
\[ \diagram 
0 \rto & \ker 1 \dto^2 \rto & \O_X \rto^1 \dto^3 & \O_Z \dto^4 \rto & 0 \\ 
0 \rto & f_* \J \rto & f_* \O_{\P W \times \P W} 
\rto & g_* \O_{\P W} \rto & 0 \enddiagram \] 
(Since $f$ is a finite morphism, $f_*$ is exact.) Since 
$4$ is an isomorphism, $\coker 2 = \coker 3$, giving a surjection 
$f_* \J \lra \Q$. Since $f^*$ is right exact, 
$f^* f_* \J \lra f^*\Q$ is also surjective. But then $q$ itself 
must be surjective, which is Claim 2. 

Now apply $\delta^*$ to $q$, then we get a surjection 
$\delta^* \J \lra g^* \Q$. By definition 
$\delta^* \J = \Omega^1_{\P W}$, hence we have an extension 
\begin{equation} 0 \lra \A \lra \Omega^1_{\P^n} \lra g^* \Q \lra 0, 
\label{ker=A} \end{equation} 
for some $\O_{\P^n}$-module $\A$. 

Let $r$ denote an integer. Tensor 
(\ref{ker=A}) by $\O_{\P^n}(r)$, and pass to the long exact sequence in 
cohomology. Assume $r \gg 0$, so that $H^1(\P^n,\A(r))=0$. 
By the Borel-Weil-Bott theorem (see \cite[p.~687]{Porras})
$H^0(\P^n,\Omega^1(r))$ is an irreducible $SL(V)$-module. 
It surjects onto $H^0(\P^n,g^*\Q(r))$, then Schur's lemma
implies that the kernel of this surjection is zero. Thus 
$H^0(\P^n,\A(r))=0$ for $r \gg 0$, which forces 
$\A = 0$. The proposition is proved. \qed 

\medskip 

\begin{Lemma} \sl 
Let $r \in \Z$. Then the group 
$H^q(\P^N, \Q(r))$ is nonzero for at most one 
value of $q$. Specifically, the only such cases are the following: 
\[ \begin{array}{lll}
H^0 & = S_{(rd-1,1)} \, V & \text{for $rd \ge 2$,} \\
H^1 & = \complex     & \text{for $r =0$,} \\
H^n & = S_{(1-rd-n,1,\dots,1,0)} \, W & \text{for $rd \le -n$.} 
\end{array} \]
\label{lemma.omegaH} \end{Lemma}
\noindent Here $S_\lambda(-)$ is the Schur functor associated to the 
partition $\lambda$ (see \cite[Ch.~6]{FH}). 

\smallskip 
\demo 
We have an isomorphism $g^* \, \Q(r) = \Omega^1_{\P W} \otimes 
\O_{\P W}(rd)$, and then the result follows from the
Borel-Weil-Bott theorem. \qed 
\begin{Lemma} \sl 
For $q \ge 1$, the group 
\[ H^q(\P^N, f_* \O_{\P W \times \P W}(r)) \] 
is nonzero iff $q = 2n$ and $r < -\frac{n}{e}$. 
\label{lemma.f*H} \end{Lemma}
\demo 
From the Leray spectral sequence and the 
K{\"u}nneth formula, 
\[ \begin{aligned} 
{} & H^q(\P^N, f_* \O_{\P W \times \P W}(r)) = \\ 
\bigoplus\limits_{i+j=q} & 
H^i(\P W, \O_{\P^n}(rd-re)) \otimes 
H^j(\P W, \O_{\P^n}(re)). 
\end{aligned} \] 
The summand $H^i \otimes H^j$ is nonzero, iff $i=j=n$ and 
the twist in each factor is $ < -n$ (see \cite[Ch.~III.5]{Ha}). 
Since $e < d-e$, the claim follows. \qed

\section{The regularity of $X$} 
Now we come to the proof of the Theorem \ref{main.theorem}. 
Henceforth we always assume that $m,q$ are positive integers in the 
range 
\begin{equation}  m \ge 0, \quad 1 \le q \le N. 
\label{mq.range} \end{equation} 
Define the predicate 
\[ \R(m,q): \; H^q(\P^N, \I_X(m-q)) = 0. \] 
Consider the following three conditions on $m$: 
\begin{itemize} 
\item[C1.] $m \ge 4$, 
\item[C2.] $d \, (m-n-2) \ge 1-n$, 
\item[C3.] $m-2n-1 \ge -\frac{n}{e}$. 
\end{itemize} 
We will reformulate the main theorem as follows: 
\begin{Theorem} \sl 
Let $m$ be a fixed positive integer. 
Then $\R(m,q)$ is true iff C1--C3 are satisfied. 
\end{Theorem} 
The smallest integer $m_0$ satisfying C1--C3 is the one defined by 
formula (\ref{defn.m0}) from the introduction. 

We will use a spectral 
sequence argument which will prove the theorem for all $q > 2$. 
This part merely amounts to checking that C1--C3 annihilate the 
$E^1$ terms in the correct positions, 
which necessarily makes somewhat tedious reading. The cases 
$q=1,2$ will follow from Propositions~\ref{q=2}~and~\ref{q=1} below, and 
the proof of the latter proposition will be completed only in 
Section~\ref{section.proof.I.II}. 

\medskip 

\demo 
We splice the exact sequence (\ref{conductor.seq}) with 
\[ 0 \lra \I_X \lra \O_{\P^N} \lra \O_X \lra 0,  \] 
and get a complex 
\[ \C^\bullet: 0 \ra \C^0 \ra \C^1 \ra \C^2 \ra 0, \] 
where 
$\C^0 = \O_{\P^N}, \C^1 = f_*\O_{\P^n \times \P^n}, \C^2 = \Omega^1_Z$. 
By construction, $H^0(\C^\bullet) = \I_X$, and 
$H^a(\C^\bullet) = 0$ for $a=1,2$. We have a spectral sequence 
\[ \begin{aligned} E_1^{a,b} & = H^b(\C^a(m-q)), \quad 
\quad \delta_r^{a,b} = E_r^{a,b} \lra E_r^{a+r,b-r+1} \\ 
E_\infty^{a,b} & \Rightarrow H^{a+b}(\I_X(m-q)); 
\end{aligned} \] 
in the range $0 \le a \le 2, 0 \le b \le N$. We will refer to this 
spectral sequence as $\Sigma_{m-q}$. 

Given the conditions (\ref{mq.range}), 
Lemmata~\ref{lemma.omegaH} and~\ref{lemma.f*H} imply that 
all the entries in $E_1$  away from the points 
\[ (a,b) = (0,0),(1,0),(2,0),(1,2n),(2,1),(2,n) \] 
are zero. This forces $E_2 = E_\infty$. To check the truth of 
$\R(m,q)$, we look at the terms $E_1^{a,b}$ in $\Sigma_{m,q}$ which are 
on the line $a +b=q$. Thus, for 
\begin{equation} q \notin \{1,\, 2,\, 3,\, n+2,\, 2n+1\}, 
\label{list.q} \end{equation}
all the terms on this line are zero, and hence $\R(m,q)$ is true. 

Firstly assume $n >1$, then the numbers in~(\ref{list.q}) are 
all distinct. Now $\R(m,n+2)$ holds iff $E_1^{2,n}=0$, 
and the latter is equivalent to C2 by Lemma~\ref{lemma.omegaH}. 
Similarly, $\R(m,2n+1) \iff E_1^{1,2n}=0 \iff$C3 by 
Lemma \ref{lemma.f*H}. Now C3 implies 
$m \ge 3$, and then $\R(m,3) \iff E_1^{2,1}=0 \iff m \neq 3$. We have 
shown that $\R(m,q)$ is true for $q \ge 3$ iff C1--C3 hold. Hence it is 
enough to show that $\R(m,2),\R(m,1)$ hold for $m \ge 4$. These 
claims will follow from Propositions \ref{q=2} and \ref{q=1}  
respectively. 

Now assume $n=1$, then the list in~(\ref{list.q}) is $\{1,2,3\}$. 
Assume $q=3$, and consider the entries 
$E_1^{a,b}$ for $(a,b)=(0,3),(1,2),(2,1)$. The first is always zero (since 
$m \ge 0$), and the rest are zero iff $m \ge 4$. Thus 
C1 holds (which entails C2,C3) iff $\R(m,3)$ holds. 
This leaves us with $q=1,2$, and again we are done by 
Propositions \ref{q=2} and \ref{q=1}. \qed 

\begin{Proposition} \sl 
Let $r \ge 1$. Then the morphism 
\[ \alpha_r: H^0(f_*\O_{\P^n \times \P^n}(r)) \lra 
 H^0(\Q(r)) \] 
is surjective. \label{q=2} 
\end{Proposition}
Since $\alpha_{m-2}$ is the morphism $\delta_1^{1,0}$ of $\Sigma_{m,2}$, 
its surjectivity implies that $E_2^{2,0}=0$, i.e., $\R(m,2)$ 
holds for $m \ge 3$. 

\smallskip 

\demo 
At the level of representations, the morphism is 
\[ \alpha_r: S_{r(d-e)} \otimes S_{re} \lra S_{(rd-1,1)}. \] 
The target of $\alpha_r$ is an irreducible $SL(V)$-module, hence  $\alpha_r$ is either surjective or zero by 
Schur's lemma. 
For $r \ge 1$, the sheaf $\Q(r) = \Omega^1_{\P^n}(rd)$ is generated by 
global sections, hence the latter is impossible. This shows that 
$\alpha_r$ is surjective. \qed 

\medskip 

Finally, consider the spectral sequence $\Sigma_{m,1}$ with $m \ge 4$. 
The truth of $\R(m,1)$ will follow if we can show that $\delta_1^{1,0}$ surjects onto the kernel 
of $\delta_1^{2,0}$. This is the content of the following proposition: 
\begin{Proposition} \sl 
Let $r \ge 3$. Then the morphism 
\[ \beta_r: H^0(\O_{\P^N}(r)) \lra H^0(\O_X(r)) \] 
is surjective. 
\label{q=1} \end{Proposition}
\demo 
For ease of reference, let us define the set 
\begin{equation} 
\A_r = \{p: 0 \le p \le re, \, p \neq 1\}. 
\label{Ar} \end{equation} 
Now the target of $\beta_r$ is 
\begin{equation} H^0(\O_X(r)) = \ker \alpha_r = 
\bigoplus\limits_{p \in \A_r} S_{(rd-p,p)}. 
\label{oxr} \end{equation}
(This follows from the Littlewood-Richardson 
rule, see~\cite[Appendix A]{FH}.) 
Let $\pi_p$ denote the projection onto the 
$p$-th summand. Then, $\pi_p \circ \beta_r$ is equal to the composite 
\begin{equation}
\begin{aligned} 
S_r(S_d) \stackrel{1}{\lra} S_r(S_{d-e} \otimes S_e) & 
\stackrel{2}{\lra} S_r(S_{d-e}) \otimes S_r(S_e) \\ 
& \stackrel{3}{\lra} S_{r(d-e)} \otimes S_{re} 
\stackrel{\pi_p}{\lra} S_{(rd-p,p)}. 
\end{aligned} \label{beta.r} \end{equation}
The map 1 is given by applying $S_r(-)$ to the coproduct 
map, 2 comes from the `Cauchy decomposition' (see \cite{ABW}), 
and 3 is the `multiplication' map. 
We will show that for $p \in \A_r$, the 
map $\pi_p \circ \beta_r$ is not identically zero, and hence 
surjective. Since the cokernel of $\beta_r$ is a direct summand of 
the target of $\beta_r$, this will prove that the cokernel is zero. 
By the argument of \cite[p.~11]{AC1}, it is enough to show 
this when $\dim V = 2$. We defer the proof to the next 
section. \qed 

\smallskip 

As a corollary to the proposition, we get a formula 
for the character of the degree $r$ part of $I_X$. 
\begin{Corollary} \sl 
For  $r \ge 3$, we have an equality 
\[ [(I_X)_r] = 
[S_r(S_d)]- \sum\limits_{p \in \A_r} \, [S_{(rd-p,p)}],  \] 
where $[-]$ denotes the formal character of an $SL(V)$-representation. 
\label{gr} \qed \end{Corollary} 
For $r=2$, we have the formula 
$[(I_X)_2] = \sum\limits_{e < p \le \frac{d}{2}} [S_{(2d-2p,2p)}]$. 
Hence the ideal has quadratic generators except when $d$ is odd and 
$e = \frac{d-1}{2}$.

\section{Transvectants} \label{section.trans}
We resume the proof of Proposition \ref{q=1}, under the hypothesis 
$\dim V = 2$. 
The argument is by induction on $r$, and the passage from 
$r$ to $r+1$ uses the symbolic calculus on binary forms. 
The notations will be consistent with 
those of Section~\ref{section.symbolic.method}. 

Given binary forms $A,B$ of degrees $a,b$ in variables 
$\ux = (x_0,x_1)$, their $k$-th transvectant $(A,B)_k$ is 
defined by formula~(\ref{transvectdef}). 
It is the image of $A \otimes B$ via the projection 
\[ S_a \otimes S_b \lra S_{a+b-2k}. \] 

\subsection{} \label{section.defn.theta}
Now define the predicate 
\begin{equation} \Theta(r,p): \pi_p \circ \beta_r \neq 0. \end{equation}
We want to show that $\Theta(r,p)$ holds for all $r \ge 3, p \in \A_r$. 
Consider the commutative diagram 
\[ \diagram 
S_r(S_d) \otimes S_d \dto \rto & S_{(rd-p,p)} \otimes S_d 
\dto^{u_r^{(p,p')}} \\ 
S_{r+1}(S_d)  \rto & S_{(rd+d-p',p')} \enddiagram \] 
where the horizontal maps are $(\pi_p \circ \beta_r) \otimes \text{id}$, and 
$\pi_{p'} \circ \beta_{r+1}$ respectively, and the vertical 
map $u_r^{(p,p')}$ is the composite 
\[ S_{(rd-p,p)} \otimes S_d \lra H^0(\O_X(r)) \otimes S_d \lra 
H^0(\O_X(r+1)) \lra S_{(rd+d-p',p')}. \]

Now consider the following two statements: 
\begin{itemize} 
\item[I.] For any $p' \in \A_3$, there exists an {\sl even} 
integer $p \in \A_2$ such that $u_2^{(p,p')} \neq 0$. 
\item[II.] Assume $r\ge 3$. Then for any $p' \in \A_{r+1}$, 
there exists a $p\in \A_r$ such that $u_r^{(p,p')} \neq 0$. 
\end{itemize} 

We claim that (I) and (II) imply 
$\Theta(r,p)$ for all $r \ge 3, p \in \A_r$. 
By~\cite[Proposition~6.2]{AC1}, 
$\Theta(2,p)$ holds for all even $p \in \A_2$. Now assume the result for 
$r$, and let $p' \in \A_{r+1}$. Let $p$ be an integer whose 
existence is guaranteed by either (I) or (II), depending on 
whether $r$ equals or exceeds $2$. By hypothesis 
$\pi_p \circ \beta_r$ is surjective, hence the composite 
\[ u_r^{(p,p')} \circ \{ (\pi_p \circ \beta_r) \otimes \text{id}\} \] 
is nonzero; this forces $\pi_{p'} \circ \beta_{r+1} \neq 0$. 

\subsection{} 
It remains to prove (I) and (II). The map $u_r^{(p,p')}$ is defined as the 
composite 
\begin{equation} \begin{aligned} 
{} & S_{rd-2p} \otimes S_d \stackrel{1}{\lra} 
(S_{r(d-e)} \otimes S_{re}) \otimes (S_{d-e} \otimes S_e) \stackrel{2}{\lra} \\ 
& S_{(r+1)(d-e)} \otimes S_{(r+1)e} \stackrel{3}{\lra} S_{(r+1)d-2p'} \, , 
\end{aligned} \label{ur} \end{equation} 
where $1$ is the tensor product of two coproduct maps, 
$2$ is obtained by regrouping, and $3$ is the projection. 

Now let $A,B$ denote binary forms of degrees $rd-2p,d$ respectively. 
We will now follow the component maps in (\ref{ur}), and get a 
step-by-step procedure
for calculating the image $u_r^{(p,p')}(A \otimes B)$. 
Introduce new variables $\uy =(y_0,y_1)$, and 
let 
\[ \Lambda = 
\sum\limits_{i=0}^{r(d-e)} \binom{r(d-e)}{i} \, 
l_i \, x_0^{r(d-e)-i}x_1^i, \quad M = 
\sum\limits_{j=0}^{re} \binom{re}{j} \, m_j \, x_0^{re-j}x_1^j, 
\] 
denote {\sl generic} binary forms of degrees $r(d-e),re$. 
(That is to say, the $l,m$ are thought of as independent indeterminates.)
\begin{itemize} 
\item 
Let $T_1 = (\Lambda,M)_p$, and $T_2 = (A,T_1)_{rd-2p}$. 
Then $T_2$ does not involve $x_0,x_1$. 
\item 
Obtain $T_3$ by making the substitutions 
\[ l_i = x_1^{r(d-e)-i}(-x_0)^i, \quad 
   m_j = y_1^{re-j}(-y_0)^j \] 
in $T_2$. 
\item 
Let 
\[ T_4 = (y_0 \, \frac{\partial}{\partial x_0} + y_1 \, 
\frac{\partial}{\partial x_1})^e \, B, \] 
usually called a partial polarization of $B$. By construction, 
$T_3$ and $T_4$ have respective bidegrees $(rd-re,re)$ 
and $(d-e,e)$ in the sets $\ux,\uy$. 
\item 
Let $T_5 = T_3 \, T_4$, and $T_6 = \Omega_{\ux\uy}^{p'} \, T_5$. 
\item 
Finally $u_r^{(p,p')}(A \otimes B)$ is obtained by 
substituting $x_0,x_1$ for $y_0,y_1$ in $T_6$. 
\end{itemize} 

\subsection{} A translation of this construction
into the classical symbolic calculus, according to 
Section~\ref{section.symbolic.method}, 
amounts to the following algebraic calculations with 
differential operators, to be understood {\sl verbatim}. 
Write
\[
\begin{array}{ll}
A(\ux)=\int_A {\rm d}a\ a_{\ux}^{rd-2p}, &
B(\ux)=\int_B {\rm d}b\ b_{\ux}^{d}, \\
\Lambda(\ux)=\int_\Lambda {\rm d}\lambda \, 
\lambda_{\ux}^{r(d-e)}, &
M(\ux)= \int_M {\rm d}\mu\ \, \mu_{\ux}^{re}. \end{array}
\]
Now, from the discussion in Section~\ref{section.symbolic.method},
\[
T_1=\int_\Lambda {\rm d}\lambda \int_M {\rm d}\mu \, 
(\lambda \, \mu)^p\ \lambda_\ux^{r(d-e)-p} \, 
\mu_\ux^{re-p}, \]
and therefore
\[ \begin{aligned} 
T_2 & =\frac{1}{(rd-2p)!^2} 
\left. \left\{ 
\Omega_{\ux\uy}^{rd-2p} \, 
A(\ux) \, T_1(\uy) \right\}
\right|_{\uy:=\ux}, \\ 
& = \frac{1}{(rd-2p)!^2}
\left. \left\{ 
\Omega_{\ux\uy}^{rd-2p} \, 
\int_A {\rm d}a \int_\Lambda {\rm d}\lambda 
\int_M {\rm d}\mu \ (\lambda \, \mu)^p \, 
a_\ux^{rd-2p} \, \lambda_\uy^{r(d-e)-p} \, 
\mu_\uy^{re-p} \right\}
\right|_{\uy:=\ux} \\ 
& = \frac{1}{(rd-2p)!^2} 
\int_A {\rm d}a 
\int_\Lambda {\rm d}\lambda
\int_M {\rm d}\mu \ (\lambda \, \mu)^p 
\left. \left\{ J(\ux,\uy)
\right\} \right|_{\uy:=\ux}, 
\end{aligned} \] 
where
\[ \begin{aligned} 
J(\ux,\uy) & =\Omega_{\ux\uy}^{rd-2p} \, 
a_\ux^{rd-2p} \, \lambda_\uy^{r(d-e)-p} \, \mu_\uy^{re-p} \\  
& =(rd-2p)! \lp a_0\frac{\partial}{\partial y_1}
-a_1\frac{\partial}{\partial y_0} \rp^{rd-2p} 
\lambda_\uy^{r(d-e)-p} \, \mu_\uy^{re-p} \\ 
& =(rd-2p)!^2 \, (a \, \lambda)^{r(d-e)-p} \, (a \, \mu)^{re-p}, 
\end{aligned} \]
which (as expected) does not involve $\ux$ or $\uy$. Therefore
\[
T_2=\int_A {\rm d}a \int_\Lambda {\rm d}\lambda
\int_M {\rm d}\mu \, (\lambda \, \mu)^p \, 
(a \, \lambda)^{r(d-e)-p} \, (a \, \mu)^{re-p}. 
\]
Now the substitution $l_i = x_1^{r(d-e)-i}(-x_0)^i$
implies that the differential operator $\int_\Lambda {\rm d}\lambda$
can be rewritten as
\[ \begin{aligned}
\int_\Lambda {\rm d}\lambda & =
\frac{1}{[r(d-e)]!} \sum_{i=0}^{r(d-e)}
\binom{r(d-e)}{i} x_1^{r(d-e)-i} (-x_0)^i \ \frac{\partial^{r(d-e)}}
{\partial \lambda_0^{r(d-e)-i} \, \partial \lambda_1^{i}} \\ 
& = \frac{1}{[r(d-e)]!} 
\lp x_1\frac{\partial}{\partial \lambda_0} 
- x_0\frac{\partial}{\partial \lambda_1} \rp^{r(d-e)} . 
\end{aligned} \] 
Likewise, from the substitution $m_j = y_1^{re-j}(-y_0)^j$, we get 
\[ \int_M {\rm d}\mu= \frac{1}{(re)!}
\lp y_1\frac{\partial}{\partial \mu_0} - y_0\frac{\partial}{\partial \mu_1}
\rp^{re}. \]
Now an easy calculation gives
\[
\begin{aligned}
{} & \int_M {\rm d}\mu \, (\lambda \, \mu)^p (a \, \mu)^{re-p} \\
& =\frac{1}{(re)!} \lp y_1\frac{\partial}{\partial \mu_0}
- y_0\frac{\partial}{\partial \mu_1} \rp^{re}
\lp \lambda_0 \, \mu_1-\lambda_1 \, \mu_0
\rp^p \lp a_0 \, \mu_1-a_1 \, \mu_0 \rp^{re-p} \\ 
& = \lp -y_1 \, \lambda_1-y_0 \, \lambda_0 \rp^p \lp
-y_1 \, a_1-y_0 \, a_0 \rp^{re-p} 
=(-1)^{re} \, \lambda_\uy^p \, a_\uy^{re-p}.  \end{aligned} \]
As a result, 
\[ \begin{aligned}
{} & \int_\Lambda {\rm d}\lambda
\int_M {\rm d}\mu \, (\lambda \, \mu)^p
(a \, \lambda)^{r(d-e)-p} \, (a \, \mu)^{re-p}\\
&  =\frac{(-1)^{re} \, a_\uy^{re-p}}{[r(d-e)]!}
\lp x_1\frac{\partial}{\partial \lambda_0} -
x_0\frac{\partial}{\partial \lambda_1} \rp^{r(d-e)}
\lp a_0\lambda_1-a_1\lambda_0 \rp^{r(d-e)-p} 
\lp \lambda_0 y_0+\lambda_1 y_1 \rp^p \\
& =(-1)^{re} \, a_\uy^{re-p} \lp -x_1 \, a_1-x_0 \, a_0 \rp^{r(d-e)-p}
\lp x_1 \, y_0-x_0 \, y_1 \rp^p, 
\end{aligned} \]
i.e., 
\[ T_3=\int_A {\rm d}a \ (-1)^{rd} \, a_\uy^{re-p} \, 
a_\ux^{r(d-e)-p} \, (\ux \, \uy)^p . \]
Now 
\[ \begin{aligned} 
T_4 & = \lp 
y_0 \frac{\partial}{\partial x_0} + y_1 \frac{\partial}{\partial x_1}
\rp^e \int_B {\rm d}b\ b_\ux^d  \\ 
& = \int_B {\rm d}b \lp y_0 \frac{\partial}{\partial x_0} 
+ y_1 \frac{\partial}{\partial x_1} \rp^e \lp
b_0 \, x_0+ b_1 \, x_1 \rp^d \\ 
& =\frac{d!}{(d-e)!} \int_B {\rm d}b \, b_\ux^{d-e} \, b_\uy^e\ ; 
\end{aligned} \] 
from which one obtains
\[ T_6= \frac{(-1)^{rd}d!}{(d-e)!} \, 
\left. \left\{ \Omega_{\ux\uy}^{p'} \int_A {\rm d}a
\int_B {\rm d}b \, (\ux \, \uy)^p \, a_\ux^{r(d-e)-p} \, 
a_\uy^{re-p} \, b_\ux^{d-e} \, b_\uy^e \right\} 
\right|_{\uy:=\ux}, \] 
or, 
\[ u_r^{(p,p')}(A \otimes B) = 
\frac{(-1)^{rd}d!}{(d-e)!} \int_A {\rm d}a 
\int_B {\rm d}b \; \E(r;p,p'), \]
where the `integrand' 
\[ \E(r;p,p')= \left.
\left\{ \Omega_{\ux\uy}^{p'} \,  (\ux\uy)^p \, 
a_\ux^{r(d-e)-p} \, b_\ux^{d-e} \, a_\uy^{re-p} \, 
b_\uy^e \right\} \right|_{\uy:=\ux} \]
is an ordinary polynomial in the variables
$a_0, a_1,b_0,b_1,x_0,x_1$. 

Note that if $\E(r;p,p')$ vanishes identically, so does
$u_r^{(p,p')}(A \otimes B)$ for any forms $A$ and $B$;
since differentiating zero gives zero. Conversely, if 
the map $u_r^{(p,p')}$ vanishes, then by applying it to 
forms $A,B$ which {\sl truly} are powers of generic linear forms, 
it would follow that $\E(r;p,p')$ itself must vanish identically. 
As a result, statements (I) and (II) from 
Section~\ref{section.defn.theta} respectively 
translate into the following: 
\begin{Proposition} \sl 
\begin{enumerate} 
\item 
For any $p' \in \A_3$, there exists an even integer 
$p \in \A_2$ such that $\E(2;p,p')$ is not identically zero. 
\item
Assume $r\ge 3$. Then for any $p' \in \A_{r+1}$, 
there exists an integer 
$p\in \A_r$ such that $\E(r;p,p')$ is not identically zero. 
\end{enumerate} \label{prop.I.II} \end{Proposition} 

The proof will be given in the next two sections. We will calculate 
the quantity $\E(r;p,p')$ explicitly and show that its 
nonvanishing is equivalent to that of a numerical 
combinatorial sum (later denoted $\cS$). 
Finally we break up the hypotheses into several subcases, and 
verify that one can always choose $p$ such that $\cS \neq 0$. 

Broadly speaking, what one has to show is that a functorially 
defined construction in multilinear algebra gives a nontrivial 
result. A thematically similar idea involving Laguerre polynomials 
appears in \cite[p.~57ff]{Shepherd-Barron}. 

%%%%%%%%%%%%%%%%%%%%%%%%%%%%%%%%%%%%%%%%%%%

\section{Transvectants of monomials 
and the quantum theory of angular momentum} \label{section.monomialtrans}
Let 
\[ L_1(\ux)=a_0 \, x_0+a_1 \, x_1, \quad L_2(\ux)=b_0 \, x_0+b_1 \, x_1, \] 
be two generic binary linear forms.
Let $\al_1,\al_2,\beta_1,\beta_2$, and $k$ be nonnegative integers such
that
$k\le \min\{\al_1+\al_2,\beta_1+\beta_2\}$.
The object of this section is to give a formula for the transvectant
\[ \T=(L_1^{\alpha_1} \, L_2^{\alpha_2}, 
L_1^{\beta_1} L_2^{\beta_2})_k \] 
which will be used later, and also to clarify its 
connection with the quantum theory of angular momentum
(see~\cite{Biedenharn}) which was alluded to
in~\cite{AC1}.

Since $\T$ is a joint covariant of $L_1(\ux)=a_\ux$ and
$L_2(\ux)=b_\ux$, it is a linear combination of bracket monomials
$(a \, b)^{i_1} \, a_\ux^{i_2} \, b_\ux^{i_3}$.
By a simple degree count,  we have 
\[ \begin{array}{ll}
i_1+i_2=\al_1+\beta_1, & i_1+i_3=\al_2+\beta_2, \\ 
i_2+i_3=\al_1+\beta_1+\al_2+\beta_2-2k. 
\end{array} \] 
This implies $i_1=k$, $i_2=\al_1+\beta_1-k$, $i_3=\al_2+\beta_2-k$.
Therefore 
\begin{equation} \T= \cN\left[
\begin{array}{c}
\alpha_1,\alpha_2\\ \beta_1,\beta_2\\ k
\end{array}
\right] (a \, b)^k \, a_\ux^{\al_1+\beta_1-k} \, 
b_\ux^{\al_2+\beta_2-k}, \label{expr.T} \end{equation} 
where 
$\cN\left[ \begin{array}{c}
\alpha_1,\alpha_2\\ \beta_1,\beta_2\\ k \end{array} \right]$ 
is a purely numerical quantity. Now specialize to $a_0=1$,
$a_1=0$, $b_0=0$, and $b_1=1$, when the right hand side becomes
\[
\cN\left[ \begin{array}{c} 
\al_1,\al_2\\ \beta_1,\beta_2\\ k \end{array} \right]
x_0^{\al_1+\beta_1-k}x_1^{\al_2+\beta_2-k} \, . \]
By definition, 
\[\begin{aligned} 
\T & = \frac{(\al_1+\al_2-k)! \, (\beta_1+\beta_2-k) !}
{(\al_1+\al_2)! \, (\beta_1+\beta_2)!}
\left.
(\frac{\partial^2 }{\partial x_0 \, \partial y_1}
-\frac{\partial^2}{\partial x_1 \, \partial y_0})^k \, 
x_0^{\alpha_1} \, x_1^{\alpha_2} \, y_0^{\beta_1} \, y_1^{\beta_2}
\right|_{\uy:=\ux} \\ 
& =\frac{(\alpha_1+\alpha_2-k)!(\beta_1+\beta_2-k)!}
{(\alpha_1+\alpha_2)!(\beta_1+\beta_2)!} \, \times \\ 
& \sum_{i=0}^k (-1)^{k-i} \binom{k}{i} 
\left. (\frac{\partial^2}{\partial x_0 \, \partial y_1})^i 
(\frac{\partial^2}{\partial x_1 \, \partial y_0})^{k-i} \, 
x_0^{\alpha_1} \, x_1^{\alpha_2} \, y_0^{\beta_1} \, y_1^{\beta_2}
\right|_{\uy:=\ux \, .} \end{aligned} \] 
This implies that 
\[ \cN\left[
\begin{array}{c}
\alpha_1,\alpha_2 \\ \beta_1,\beta_2\\ k
\end{array} \right] =
\cS\left[ \begin{array}{c} \alpha_1,\alpha_2\\ \beta_1,\beta_2\\ k
\end{array} \right]
\times \frac{(-1)^k (\al_1+\al_2-k)!(\beta_1+\beta_2-k)!k!{\al_1}!
{\al_2}!{\beta_1}!{\beta_2}!}
{(\al_1+\al_2)!(\beta_1+\beta_2)!}
\]
where, by definition 
\begin{equation} 
\begin{aligned} {} & \cS\left[ \begin{array}{c} 
\alpha_1,\alpha_2 \\\beta_1,\beta_2 \\ k \end{array}
\right] = \\ 
& \sum_{i=\max\{0,k-\al_2,k-\beta_1\}}^{\min\{
k,\al_1,\beta_2\}} 
\frac{(-1)^i}
{i!(k-i)!(\al_1-i)!(\beta_2-i)!(\al_2-k+i)!(\beta_1-k+i)!}
\label{sumformula} \end{aligned} \end{equation} 
\subsection{} 
We now connect this to the prevalent formalism in physics. 
From the original data $\alpha_1,\alpha_2,\beta_1,\beta_2,k$, 
we define numbers $j_1, j_2, j, m_1, m_2, m$ via the relations 
\[ \begin{array}{rrl} 
j_1=\frac{1}{2}(\alpha_1 + \alpha_2), & 
j_2=\frac{1}{2}(\beta_1+\beta_2), & 
j=\frac{1}{2} (\alpha_1+\alpha_2+\beta_1+\beta_2)-k \\ 
m_1=\frac{1}{2}(\alpha_2-\alpha_1), 
& m_2=\frac{1}{2}(\beta_2-\beta_1), & 
m=\frac{1}{2}(\alpha_2-\alpha_1+\beta_2-\beta_1).
\end{array} \]
Now the so-called vector-coupling or Clebsch-Gordan coefficients
are the quantities
\[ \begin{aligned} 
{} & C_{m_1,m_2,m}^{j_1,j_2,j}=
\left[\frac{(2j+1)(j_1+j_2-j)!(j_1+j-j_2)!(j_2+j-j_1)!}
{(j_1+j_2+j+1)!} \right]^{\frac{1}{2}} \\ 
& \times \left[
(j_1-m_1)!(j_1+m_1)!(j_2-m_2)!(j_2+m_2)!(j-m)!(j+m)!
\right]^{\frac{1}{2}} \times
\sum \frac{(-1)^i}{z}, \end{aligned}  \] 
where $z$ stands for 
\[ {i!(j_1+j_2-j-i)!(j_1-m_1-i)!(j_2+m_2-i)!
(j-j_2+m_1+i)!(j-j_1-m_2+i)!},\] 
and the last summation is quantified over 
\[ \min\{j_1+j_2-j,j_1-m_1,j_2+m_2\} \le i 
\le \max\{0,j_2-j-m_1, j_1-j+m_2\}. \] 
Now, using our original data, 
\[ \begin{aligned} {} & C_{m_1,m_2,m}^{j_1,j_2,j}= 
\cN\left[ \begin{array}{c}
\alpha_1,\alpha_2 \\ \beta_1,\beta_2\\ k \end{array} \right] \\ 
& \times (-1)^k\sqrt{\alpha_1+\alpha_2+\beta_1+\beta_2-2k+1} \, 
(\alpha_1+\alpha_2)! \, (\beta_1+\beta_2)! \\ 
& \times \left[
\frac{(\al_1+\beta_1-k)!(\al_2+\beta_2-k)!}
{(\al_1+\al_2+\beta_1+\beta_2-k+1)!
(\al_1+\al_2-k)!(\beta_1+\beta_2-k)!
k!\al_1!\al_2!\beta_1!\beta_2!}
\right]^{\frac{1}{2}} \, .
\end{aligned} \]

\smallskip 

Physicists also use related quantities called
{\sl Wigner's $3j$-symbols}
given by
\[
\lp
\begin{array}{ccc}
j_1 & j_2 & j\\ m_1 & m_2 & -m \end{array}
\rp=
(-1)^{j_1-j_2+m} \, (2j+1)^{-\frac{1}{2}} \, 
C_{m_1,m_2,m}^{j_1,j_2,j} . 
\] 
\subsection{} \label{uv}
We will record another calculation which will be useful in the 
next section. Introduce auxiliary variables 
$u_0$, $u_1$, $v_0$, $v_1$, and let 
$\{a \, \partial_u\}$ stand for $a_0 \, \frac{\partial}{\partial u_0}
+a_1 \, \frac{\partial}{\partial u_1}$ etc. Define 
\be
\begin{aligned} {} & \cU\left[ \begin{array}{c}
\alpha_1,\alpha_2 \\ \beta_1,\beta_2\\ k
\end{array} \right] (a,b,\ux) \\ 
= & \{a \, \partial_u\}^{\alpha_1} \, 
\{b \, \partial_u\}^{\alpha_2} \, \{a \, \partial_v\}^{\beta_1} \, 
\{b \, \partial_v\}^{\beta_2} \, (u \, v)^k \, 
u_\ux^{\al_1+\al_2-k} \, v_\ux^{\beta_1+\beta_2-k} \, .
\end{aligned}
\label{udef}
\ee 
The following result will be needed.

\begin{Lemma} \sl
\[ \cU\left[ \begin{array}{c}
\alpha_1,\alpha_2\\ \beta_1,\beta_2\\ k
\end{array} \right] (a,b,\ux)=
(\al_1+\al_2)!(\beta_1+\beta_2)!
\ \T
\]
where $\T$ denotes the monomial transvectant
$(a_\ux^{\alpha_1} \, b_\ux^{\alpha_2},
a_\ux^{\beta_1} \, b_\ux^{\beta_2})_k$.
\end{Lemma}

\demo
Using a graphical notation for symmetrizers as in~\cite{Cvitanovic}
would make the truth of the lemma visually obvious.
Alternatively, one can do the following.
Write each factor in (\ref{udef}) as a sum over indices with values
in $\{0,1\}$ :
\[
\begin{array}{ll}
\{a \, \partial_u\}=\sum_i a_i \, \frac{\partial}{\partial u_i},  &
(u \, v)=\sum_{i,j} u_i \, \epsilon_{ij} \, v_j,\\
u_\ux=\sum_i u_i \, x_i \;\; \text{etc.,} \end{array}
\]
where $\epsilon=(\epsilon_{ij})$ is the antisymmetric $2\times 2$ matrix
with $\epsilon_{01}=1$.
One has to use disjoint sets of indices for each individual factor in
(\ref{udef}), i.e., a total of $2p$ indices, with
$p=\al_1+\al_2+\beta_1+\beta_2$. Now expand $\cU$
completely, which gives an expression of the form
\[ \cU=\sum_{I,J} \, A_{I,J}
\ \frac{\partial}{\partial z_{i_1}}
\cdots \frac{\partial}{\partial z_{i_p}} 
z_{j_1}\ldots z_{j_p} \]
where $I=(i_1,\ldots,i_p)$ and $J=(j_1,\ldots,j_p)$
are collections of indices running from $1$ to $4$,
by definition 
\[ \uz=(z_1,z_2,z_3,z_4)=(u_0,u_1,v_0,v_1), \]
and finally $A_{I,J}$ are coefficients depending on $a,b,\ux$,
the detailed expression of which we spare the reader. Now
\[ \frac{\partial}{\partial z_{i_1}} \cdots
\frac{\partial}{\partial z_{i_p}} \, z_{j_1}\ldots z_{j_p}
= \sum_\sigma \, \prod_{\nu=1}^p \, 
\delta_{i_{\nu} j_{\sigma(\nu)}}
\]
where $\sigma$ denotes a permutation of the set $\{1,\ldots,p\}$.
Since
\[ \frac{\partial}{\partial z_{i_1}}
\cdots \frac{\partial}{\partial z_{i_p}} z_{j_1}\ldots z_{j_p}=
\frac{\partial}{\partial z_{j_1}} \cdots
\frac{\partial}{\partial z_{j_p}} z_{i_1}\ldots z_{i_p}, 
\]
and one can exchange the role of the dummy summation indices $I$ and $J$,
we get 
\[ \cU=\sum_{I,J} A_{J,I} \frac{\partial}{\partial z_{i_1}} 
\cdots \frac{\partial}{\partial z_{i_p}} z_{j_1}\ldots z_{j_p} \period 
\]
Now undo the previous expansion of $\cU$ to find
\[ \cU\left[ \begin{array}{c}
\alpha_1,\alpha_2\\ \beta_1,\beta_2\\ k
\end{array} \right] (a,b,\ux)=
(\partial_u)_\ux^{\al_1+\al_2-k} \, 
(\partial_v)_\ux^{\beta_1+\beta_2-k} \; 
\cW(u,v), \]
where 
$(\partial_u)_\ux=x_0 \, \frac{\partial}{\partial u_0}+
x_1 \, \frac{\partial}{\partial u_1}$ etc., and 
\[ 
\cW(u,v) = (\frac{\partial^2 }{\partial u_0 \, \partial v_1} - 
\frac{\partial^2}{\partial u_1 \, \partial v_0})^k \, 
a_u^{\alpha_1} \, b_u^{\alpha_2} \, 
a_v^{\beta_1} \, b_v^{\beta_2}. \] 
Since $\cW(u,v)$ is homogeneous in $u, v$ of respective degrees
$\alpha_1+\alpha_2-k$ and $\beta_1+\beta_2-k$, 
\[ \cU\left[ \begin{array}{c} \alpha_1,\alpha_2\\ 
\beta_1,\beta_2\\ k \end{array}
\right] (a,b,\ux)=
(\al_1+\al_2-k)! \, (\beta_1+\beta_2-k)! \, 
\cW(u,v)|_{u,v:=\ux}, \]
and the lemma follows. \qed

\begin{Corollary}
\[ \begin{aligned} {} & \cU\left[ \begin{array}{c}
\alpha_1,\alpha_2\\ \beta_1,\beta_2\\ k
\end{array} \right] (a,b,\ux)  \\ 
= \, & \cN\left[\begin{array}{c}
\alpha_1,\alpha_2\\ \beta_1,\beta_2\\ k \end{array} \right] \, 
(\alpha_1+\alpha_2)! \, (\beta_1+\beta_2)! \times 
(a \, b)^k \, a_\ux^{\alpha_1+\beta_1-k} \, b_\ux^{\alpha_2+\beta_2-k} \, .
\end{aligned} \] 
\end{Corollary}

%%%%%%%%%%%%%%%%%%%%%%%%%%%%%

\section{Proof of proposition \ref{prop.I.II}} 
\label{section.proof.I.II}

In this section we will prove the remaining proposition, and 
hence complete the proof of the main theorem. 
\medskip 
 
Let us write 
\[ \M=\Omega_{\ux\uy}^{p'} \, (\ux \, \uy)^p \, a_{\ux}^{r(d-e)-p} \, 
b_{\ux}^{d-e} \, a_{\uy}^{re-p} \, b_{\uy}^e,  \]
then 
\[ \E = \E(r;p,p') = \left. \M \right|_{\uy := \ux \,} \] 
is the expression to be calculated. 
Introduce pairs of variables $u_0$, $u_1$, $v_0$, $v_1$ as in 
Section~\ref{uv}. In the notation introduced there, 
\[ 
\{a \, \partial_u\}^{r(d-e)-p} \, \{b \, \partial_u\}^{d-e} \, 
u_\ux^{(r+1)(d-e)-p} 
= ((r+1)(d-e)-p)! \, a_{\ux}^{r(d-e)-p} \, b_{\ux}^{d-e}, 
\] 
and similarly 
\[ \{a \, \partial_v\}^{re-p} \, \{b \, \partial_v\}^e \, 
v_\uy^{(r+1)e-p} = ((r+1)e-p)! \, a_{\uy}^{re-p} \; b_{\uy}^e \period
\] 
Hence 
\[ \begin{aligned} 
{} & \M = \frac{1}{[(r+1)e-p]! \, [(r+1)(d-e)-p]!} \, 
\Omega_{\ux\uy}^{p'} \\ 
& \lp (\ux \, \uy)^p {\{a \, \partial_u\}}^{r(d-e)-p} \, 
\{b \, \partial_u\}^{d-e} \, \{a \, \partial_v\}^{re-p} \, 
\{b \, \partial_v\}^e \, u_\ux^{(r+1)(d-e)-p} \, v_\uy^{(r+1)e-p} \rp. 
\end{aligned} \] 
We can commute those partial differential operators
which act on disjoint sets of variables, this gives 
\[ \M= \frac{\{a \, \partial_u\}^{r(d-e)-p} \, \{b \, \partial_u\}^{d-e} \, 
\{a \, \partial_v\}^{re-p} \, \{b \, \partial_v\}^e}
{[(r+1)e-p]! \, [(r+1)(d-e)-p]!} \times \cP, \] 
where 
\[ \cP= 
\Omega_{\ux\uy}^{p'} 
\{ (\ux \, \uy)^p \, u_\ux^{(r+1)(d-e)-p} \, 
v_\uy^{(r+1)e-p} \}. \] 
\subsection{} 
Now let $m=(r+1)(d-e)-p$, $n=(r+1)e-p$, so that $n<m$. 
(Of course, this $n$ is entirely unrelated to the one from 
Section~\ref{section.notation}. The latter plays no role in this 
calculation.) By the Clebsch-Gordan series (see \cite[Ch.~IV]{GrYo}), 
\[ u_\ux^m \, v_\uy^n =
\sum_{j=0}^{n} \frac{\binom{m}{j}\binom{n}{j}}
{\binom{m+n-j+1}{j}} \, (\ux \, \uy)^j
(u_\ux^m, v_\uy^n)_{\uy^{n-j}}^j, \] 
where, by definition,
\[ (u_\ux^m, v_\uy^n)_{\uy^{n-j}}^j 
= \frac{(m-j)!}{(m+n-2j)!}
\{\uy \, \partial_\ux\}^{n-j} \, (u \, v)^j \, 
u_\ux^{m-j}  \, v_\ux^{n-j} \period \] 
Now introduce new variables $w_0,w_1$, and 
rewrite the last expression as 
\[ (u_\ux^m, v_\uy^n)_{\uy^{n-j}}^j =
\frac{(u \, v)^j}{(m+n-2j)!} \, 
\{u \, \partial_w\}^{m-j} \, \{v \, \partial_w\}^{n-j} \, 
w_{\ux}^{m-j} w_{\uy}^{n-j} \period \] 
Then $\cP$ can be written as 
\[ \begin{aligned} 
{} & \Omega_{\ux\uy}^{p'} \; 
(\ux \, \uy)^p \, u_\ux^m \, v_\uy^n 
= \Omega_{\ux\uy}^{p'}
\left[ (\ux \, \uy)^p \,
\sum_{j=0}^{n}
\frac{\binom{m}{j}\binom{n}{j}}{\binom{m+n-j+1}{j}}
(\ux \, \uy)^j (u_\ux^m, v_\uy^n)_{\uy^{n-j}}^j \right] \\
= & \; \Omega_{\ux\uy}^{p'} 
\left[ (\ux\uy)^p \sum_{j=0}^n \, \U_j \right] 
\end{aligned} \] 
where 
\[ \U_j = \frac{\binom{m}{j}\binom{n}{j}}{\binom{m+n-j+1}{j}}
\frac{(\ux \, \uy)^j}{(m+n-2j)!} 
(u \, v)^j \, \{u \, \partial_w\}^{m-j} \, 
\{v \, \partial_w\}^{n-j} \, w_{\ux}^{m-j} \, w_{\uy}^{n-j} \period 
\] 
Again the point is that the differential operators can be commuted! 
Hence 
\[ \begin{aligned} \cP = 
\sum_{j=0}^{n} & \, \frac{m!}{j! \, (m-j)!} \frac{n!}{j! \, (n-j)!} 
\times \frac{j! \, (m+n-2j+1)!}{(m+n-j+1)!} \\ 
& \times \frac{1}{(m+n-2j)!} (u \, v)^j \, 
\{u \, \partial_w\}^{m-j} \, \{v \, \partial_w\}^{n-j} \times \Q, 
\end{aligned} \] 
where 
\[ \Q= \Omega_{\ux\uy}^{p'} 
(\ux \, \uy)^{p+j} \, w_{\ux}^{m-j} \, 
w_{\uy}^{n-j} \, . \] 
The last expression occurs frequently in classical 
invariant theory. It is calculated, for instance, in~\cite[\S 3.2.6]{Glenn}. 
If $p'>p+j$ then $\Q$ is zero, and if $p'\le p+j$ then it equals 
\[ \frac{(p+j)!}{(p+j-p')!} \times
\frac{(m+n+p-j+1)!}{(m+n+p-p'-j+1)!} \, 
(\ux \, \uy)^{p+j-p'} \, w_{\ux}^{m-j} \, w_{\uy}^{n-j} \period \] 
As a result, 
\begin{equation} \begin{aligned} \cP = 
\sum_j \; & \bbbone_{\left\{{0\le j\le n}\atop{p'\le p+j}\right\}} \, 
\frac{m! \, n! \, (m+n-2j+1)}{j! \, (m-j)! \, (n-j)! \, (m+n-j+1)!} \\ 
& \times \frac{(p+j)!}{(p+j-p')!} \times 
\frac{(m+n+p-j+1)!}{(m+n+p-p'-j+1)!}  \\ 
& \times (u \, v)^j \{u \, \partial_w\}^{m-j}
\{v \, \partial_w\}^{n-j} (\ux \, \uy)^{p+j-p'} w_{\ux}^{m-j}
w_{\uy}^{n-j} \period 
\end{aligned} \label{formula.N} \end{equation} 
Here $\bbbone_{\{ \; \}}$ denotes the {\sl characteristic function} of 
that set. 
\subsection{} 
Now recall that 
\[ \M= \frac{\{a\partial_u\}^{r(d-e)-p}\{b \, \partial_u\}^{d-e}
\{a \, \partial_v\}^{re-p} \{b \, \partial_v\}^e}
{m! \, n!} \times \cP \] 
with $m=(r+1)(d-e)-p$, and $n=(r+1)e-p$. 
The quantity we are interested in is
$\E=\M|_{\uy:=\ux}$. When we set $\uy=\ux$ in (\ref{formula.N}), 
the term corresponding to $j=p'-p$ is the only one that survives. 
Therefore $\E$ equals 
\[ \begin{aligned} 
{} & \bbbone_{\{0\le p'-p\le n\}} \times 
\{a \, \partial_u\}^{r(d-e)-p} \, 
\{b \, \partial_u\}^{d-e} \, 
\{a \, \partial_v\}^{re-p} \, 
\{b \, \partial_v\}^e \; \times \\ 
& (u \, v)^{p'-p} \, \{u \, \partial_w\}^{m-p'+p} \, 
\{v \, \partial_w\}^{n-p'+p} \, w_{\ux}^{m+n-2p'+2p} \; \times  \\ 
& \frac{p'! \, (m+n+2p-p'+1)!}
{(p'-p)!(m-p'+p)!(n-p'+p)!(m+n-p'+p+1)!(m+n-2p'+2p)!} \period
\end{aligned} \]
Now
\[ \{u \, \partial_w\}^{m-p'+p} \, 
\{v \, \partial_w\}^{n-p'+p} \, 
w_{\ux}^{m+n-2p'+2p} =
(m+n-2p'+2p) \, ! \, u_\ux^{m-p'+p} \, v_\ux^{n-p'+p} \period 
\] 
The condition $p'-p\le n$ is equivalent
to the hypothesis $p'\le (r+1)e$, and can therefore be dropped. Hence 
\[ \begin{aligned} 
\E & = \bbbone_{\{p\le p'\}} \, \times \cU \, \times \\
& \frac{p'! \, ((r+1)d-p'+1)!} 
{(p'-p)!((r+1)(d-e)-p')!((r+1)e-p')!((r+1)d-p'-p+1)!} 
\end{aligned}\] 
where 
\[ \cU = \cU\left[ \begin{array}{c} r(d-e)-p,d-e\\ re-p,e\\ p'-p
\end{array} \right] \] 
in the notation of Section~\ref{uv}. 
As a result, $\E$ is nonzero iff the sum 
\[ \cS=\cS\left[ \begin{array}{c}
r(d-e)-p,d-e\\ re-p,e\\ p'-p \end{array}
\right] \]
defined by Formula (\ref{sumformula}) is nonzero. Let 
\begin{equation} \begin{aligned} 
\Low(d,e,r,p',p)= & \max \, \{0,p'-p-(d-e),p'-re\}, \\ 
\Hi(d,e,r,p',p)= & \min \, \{p'-p,e,r(d-e)-p\}, 
\end{aligned} \label{LowHi} \end{equation} 
henceforth written as $\Low$ and $\Hi$ if no confusion is likely. 
The index of summation in the definition of $\cS$ runs from
$\Low$ to $\Hi$. 

\subsection{} In~\cite[\S 6]{AC1} we were able to produce 
closed formulae for such coefficients in analogous cases, and 
then to determine whether they were nonzero. This was due to 
the particular form of these coefficients, which allowed the 
use of some standard summation theorems for hypergeometric series. 
On the contrary, we now have five independent parameters $d,e,r,p',p$, 
and no such closed formulae seem to apply.
It does not seem either that all the situations where there exist 
such summation formulae (say for terminating hypergeometric series or 
Wigner's $3nj$-symbols) have been classified in the framework
of Wilf-Zeilberger theory (cf.~\cite{KrattenthalerR,PetkovsekWZ,WilfZ}). 
Moreover, the determination of the zeros of  $3nj$-symbols (even
in the simplest case of $3j$-symbols) is an outstanding open problem in 
the quantum theory of angular momentum 
(see~\cite[Ch.~5, Topic 10]{Biedenharn2} or~\cite{RaynalVRR}). 
It is quite intriguing that some of these zeros have been given an 
explanation involving exceptional Lie groups. The issue may well be 
related to the article by Dixmier~\cite{Dixmier}
(see also~\cite[Chap. 6]{Okubo}), where, 
among a bestiary of nonassociative algebras, the octonions 
are realised by a construction using transvectants of binary forms.

\subsection{} We can conclude the proof of Proposition \ref{prop.I.II} 
only because we have some freedom in the choice of $p$. 
We break up the allowed values of $r$ and $p'$ into several 
cases, and by analysing which entries realize the 
maximum and minimum in (\ref{LowHi}), it is always possible to 
choose the value of $p$ in such a way that the sum defining $\cS$ 
has {\sl at most two terms}. In the angular momentum parlance, these 
correspond to `stretched' $3j$-symbols. The trickiest part is 
the case $r=2$, since there we have fewer choices for $p$.

\medskip 

\noindent
{\bf Case 1 :} $r=2$, $0\le p'\le 2e$ and $p'$ even.

\noindent Choose $p=p'$. Therefore $\Low= \Hi=0$, and
\[ \cS=\frac{1}{e!(2(d-e)-p')!(2e-p')!(d-e)!} \neq 0.\] 

\bigskip 

\noindent {\bf Case 2 :} $r=2$, $0\le p'\le 2e$ and $p'$ odd.

\noindent Choose $p=p'-1$. Therefore
$\Low=0,\Hi=1$, and  
\[ \cS=\frac{(-p'+1)(d-2e)}{e!(2(d-e)-p'+1)!(2e-p'+1)!(d-e)!} \neq 0.
\]

\bigskip 

\noindent
{\bf Case 3 :} $r=2$, and  $2e< p'\le \min\{2(d-e),3e\}$.

\noindent
Choose $p=2e$. Therefore
$\Low = \Hi = p'-2e$, and 
\[ \cS=\frac{(-1)^{p'-2e}}{(p'-2e)!(3e-p')!(2(d-e)-p')!(d-e)!} \neq 0.
\] 

\bigskip 

\noindent
{\bf Case 4 :} $r=2$,  $2(d-e)< p'\le 3e$, and $p'$ even.

\noindent
Choose $p=2d-p'$. Therefore 
$\Low = \Hi = p'-2e$, and
\[ 
\cS=\frac{(-1)^{p'-2e}} 
{(p'-2e)!(3e-p')!(p'-2(d-e))!(3(d-e)-p')!} \neq 0.
\] 

\bigskip 

\noindent
{\bf Case 5 :} $r=2$, $2(d-e)< p'< 3e$, and $p'$ odd.

\noindent
Choose $p=2d-p'-1$. Therefore 
$\Low=p'-2e,\Hi=p'-2e+1$, and 
\[ \cS=\frac{(-1)^{p'-2e}(d-2e)(p'+3)}
{(p'-2e+1)!(3e-p')!(p'-2(d-e)+1)!(3(d-e)-p')!} \neq 0.
\] 

\bigskip 

\noindent
{\bf Case 6 :} $r=2$,  $p'=3e$ and $p'$ odd (i.e., $e$ is odd).

\noindent Choose $p=2d-p'-1$. Therefore
$\Low = \Hi = e$, and 
\[ \cS=\frac{-1}{e!(5e-2d+1)!(3d-6e-1)!} \neq 0.
\] 

\bigskip 

\noindent {\bf Case 7 :} $r\ge 3$,  $0\le p'\le re$ and $p'\neq 1$.

\noindent Choose $p=p'$. Therefore
$\Low = \Hi = 0$, and 
\[ \cS=\frac{1}{e!(r(d-e)-p')!(re-p')!(d-e)!}\neq 0. 
\] 

\bigskip 

\noindent
{\bf Case 8 :} $r\ge 3$, and $re< p'\le \min\{r(d-e),(r+1)e\}$.

\noindent Choose $p=re$. Therefore
$\Low = \Hi = p'-re$, and 
\[ \cS=\frac{(-1)^{p'-re}}{(p'-re)!((r+1)e-p')!(r(d-e)-p')!(d-e)!}
\neq 0.\] 

\bigskip 

\noindent
{\bf Case 9 :} $r\ge 3$, and $r(d-e)< p'\le (r+1)e$.

\noindent
Choose $p=rd-p'$. Therefore $\Low = \Hi =p'-re$, and 
\[ \cS=\frac{(-1)^{p'-re}}
{(p'-re)!((r+1)e-p')!(p'-r(d-e))!((r+1)(d-e)-p')!} \neq 0. 
\] 
The proof of Proposition \ref{prop.I.II} (and hence that of 
the main theorem) is complete. \qed 

\section{Binary forms-I}
\label{section.binary.set-th} 
In this section we will construct two covariants of 
binary $d$-ics which together define the locus $X$ set-theoretically. 

\subsection{} 
Let $F$ be a binary $d$-ic. For nonnegative integers $\mu,\nu$, 
we will use the notation 
\[ F_{x_0^\mu \, x_1^\nu} = 
\frac{\partial^{\mu+\nu} \, F}{\partial x_0^\mu \, \partial x_1^\nu}. 
\] 
Now write $F = \prod\limits_i \, l_i^{\alpha_i}$, 
where $l_i$ are pairwise non-proportional linear forms, and $\alpha_i > 0$. 
Assume furthermore that $F$ is not a power of a linear form.
Define $g_F = \gcd \, (F_{x_0},F_{x_1})$. (By our assumption, 
both $F_{x_i}$ are nonzero.)

\begin{Lemma} \sl 
With notation as above, $g_F = \prod l_i^{\alpha_i-1}$. 
\label{lemma.gcd} \end{Lemma} 
\demo Evidently $g = \prod l_i^{\alpha_i-1}$ divides both the 
$F_{x_i}$, hence write $F_{x_0} = g \, A, F_{x_1} = g \, B$. Divide 
Euler's equation $d \, F = x_0 \, F_{x_0} + x_1 \, F_{x_1}$ 
by $g$, then 
$d \, \prod l_i = x_0 \, A + x_1 \, B$. If $A,B$ have a common linear 
factor, it must be one of the $l_i$, say $l_1$. But 
\[ A = \sum\limits_i \alpha_i \frac{\partial l_i}{\partial x_0} 
(\prod\limits_{j \neq i} l_j), \] 
so $l_1|A$ implies $\frac{\partial l_1}{\partial x_0}=0$. The same 
argument on $B$ leads to $\frac{\partial l_1}{\partial x_1}=0$, 
so $l_1=0$. This is absurd, hence $A,B$ can have no common factor, 
i.e.~$g = g_F$. \qed 

\smallskip 
Now define 
\[ Y = \bigcup_{0 \le e \le [\frac{d}{2}]} X^{(d-e,e)}, \] 
the locus of binary $d$-ics with at most two distinct linear factors. 
\begin{Lemma} \sl 
For a binary $d$-ic $F(x_0,x_1)$, the following are equivalent: 
\begin{enumerate} 
\item[(i)] $F \in Y$. 
\item[(ii)] The forms 
$\U = \{x_0 \, F_{x_0}, x_0 \, F_{x_1}, x_1 \, F_{x_0}, x_1 \, F_{x_1}\}$ 
are linearly dependent. 
\end{enumerate} \end{Lemma} 
\demo 
Assume (i), then $F = x_0^{d-e} \, x_1^e$ or $x_0^d$ after a change of 
variables, and (ii) is immediate. If (ii) holds, then there exist 
linear forms $l,m$ (not both zero) such that 
$l \, F_{x_0} = m \, F_{x_1}$. But then either 
$F_{x_0},F_{x_1}$ have a common factor of degree $\ge d-2$, or one of 
them is zero. In the latter case $F$ is a power of a linear form. In 
the former case, the previous lemma implies that $F$ has at most 
two distinct linear factors. \qed 

\medskip 

Define $\D(F)$ to the Wronskian of the sequence $\U$, i.e., 
\begin{equation} \D(F) = \det 
\left| \begin{array}{rrrr} 
(x_0 \, F_{x_0})_{x_0^3} & (x_0 \, F_{x_0})_{x_0^2 \,x_1} & 
(x_0 \, F_{x_0})_{x_0 \, x_1^2} & (x_0 \, F_{x_0})_{x_1^3} \\ 
(x_0 \, F_{x_1})_{x_0^3} & (x_0 \, F_{x_1})_{x_0^2 \,x_1} & 
(x_0 \, F_{x_1})_{x_0 \, x_1^2} & (x_0 \, F_{x_1})_{x_1^3} \\ 
(x_1 \, F_{x_0})_{x_0^3} & (x_1 \, F_{x_0})_{x_0^2 \,x_1} & 
(x_1 \, F_{x_0})_{x_0 \, x_1^2} & (x_1 \, F_{x_0})_{x_1^3} \\ 
(x_1 \, F_{x_1})_{x_0^3} & (x_1 \, F_{x_1})_{x_0^2 \,x_1} & 
(x_1 \, F_{x_1})_{x_0 \, x_1^2} & (x_1 \, F_{x_1})_{x_1^3} 
\end{array} \right| \label{D(F).wronskian} 
\end{equation}
It is a covariant of $F$ of degree $4$ and order $4d-12$. 
By the previous lemma, 
\begin{equation}
 F \in Y \iff \D(F) = 0. \label{D(F)} 
\end{equation}

Recall that a binary form $A(x_0,x_1)$ is a power of a 
linear form, iff its Hessian 
\[ \He \, (A) = A_{x_0^2} \, A_{x_1^2} - (A_{x_0x_1})^2 \] 
is identically zero. Moreover, for such a form {\sl all} covariants of 
degree greater than one are identically zero. 

\subsection{} 
Now fix an integer $1 \le e \le \frac{d}{2}$, and define a {\sl rational} 
covariant 
\[ \A_e(F) = \frac{F^{2d-2e-2}}{\He(F)^{d-e}} \] 
Assume $F \in Y,\He(F) \neq 0$. By a change of variable we may write 
$F = x_0^{d-f} \, x_1^f$, for some $1 \le f \le \frac{d}{2}$. 
Then up to a nonzero multiplicative factor, 
$\He(F) = x_0^{2d-2f-2} \, x_1^{2f-2}$, and by a direct substitution 
\[ \A_e(F) = x_0^{2f-2e} \, x_1^{2d-2e-2f}. \] 
Hence $\A_e(F)$ can be a power of a linear form, iff either 
$f=e$ or $e+f=d$, i.e., iff $F \in X^{(d-e,e)}$. Hence we have proved the 
following: 
\begin{Proposition} 
\sl A binary $d$-ic $F$ (which is not a $d$-th 
power of a linear form) lies in $X^{(d-e,e)}$, iff 
$\D(F) = \He \, (\A_e(F)) = 0$. 
\end{Proposition} 

However, this criterion is not aesthetically satisfactory insomuch as 
it appeals to a rational (as opposed to a polynomial) covariant. 
To amend this, we will deduce a formula for the Hessian of a quotient 
of two forms and then apply it to $\A_e$. 

\subsection{The Hessian of a quotient} 
Let $P,Q$ denote generic binary forms of degrees $p,q \ge 0$ 
respectively. (By convention, $1$ is the generic degree zero form.) 
Define 
\[ \begin{aligned} 
z_1 & = \frac{p^2 \, (p-1) \, (2p-2q-1) \, (p-q-1)}{2(2p-1)}, \\ 
z_2 & = \frac{q^2 \, (q-1) \, (2p-2q+1) \, (p-q-1)}{2(2q-1)}, \\  
z_3 & = p \, q \, (p-q-1). \end{aligned} \] 
\begin{Theorem} \sl 
With notation as above, we have the following formal identity:
\[ \He \, (\frac{P}{Q}) = \frac{J(P,Q)}{Q^4}, \] 
where 
\begin{equation}
J(P,Q) = z_1 \, Q^2 \, (P,P)_2 + z_2 \, P^2 \, (Q,Q)_2 
+ z_3 \, (P^2,Q^2)_2. 
\label{jpq} \end{equation}
\end{Theorem} 
\demo If either $p$ or $q$ is zero, then the theorem reduces to 
an easy calculation, hence we may assume $p,q \ge 1$. 
Let $U = \frac{P}{Q}$, first
we will show that $J = Q^4 \, He(U)$ is a polynomial. 
Let us write 
\[ \partial_i = \frac{\partial}{\partial x_i}, \quad 
   \partial_{ij} = \frac{\partial^2}{\partial x_i \, \partial x_j}, \] 
then by quotient rule, 
\begin{equation} \begin{aligned} 
\partial_{ij} \, U & = 
\frac{\partial_{ij}P}{Q} - 
\frac{\partial_i P \, \partial_j Q + \partial_j P \, \partial_i Q}{Q^2} 
- \frac{P \, \partial_{ij}Q}{Q^2} + 
 \frac{2P \, \partial_iQ \, \partial_j Q}{Q^3} \\ 
& = e_1(i,j) - e_2(i,j) - e_3(i,j) + e_4(i,j). \end{aligned} 
\label{partialijf} \end{equation} 
Here $e_\star(i,j)$ are simply names for those consecutive 
expressions. Now 
\[ \He \, (U) = (\partial_{0,0} \, U) \, (\partial_{1,1} \, U)  - 
(\partial_{0,1} \, U)^2 \] 
is a linear combination of terms 
\[ E(a,b) = e_a(0,0) \, e_b(1,1) + e_b(0,0) \, e_a(1,1) 
- 2 \, e_a(0,1) \, e_b(0,1), \] 
for $1 \le a,b \le 4$. The terms $E(2,4),E(4,4)$ are zero, 
so the only term with a (possible) denominator of $Q^5$ is $E(3,4)$. 
Let us write 
\[ E(3,4) = -\frac{2P^2}{Q^5} \, E'(3,4), \] 
we will show that in fact $Q$ divides $E'(3,4)$. Now we have an 
identity 
\begin{equation} 
q^2 \, (2q-1)(q-1) \, (Q^2,Q)_2 = q^2 (q-1)^2 \, Q \, (Q,Q)_2 + 
E'(3,4); \label{q.q^2.2} \end{equation}
this follows by directly calculating the left hand side with 
formula~(\ref{trans.formula}). 
So far the entire argument works if $Q$ is any sufficiently 
differentiable function of $x_0,x_1$. But now we can 
use the homogeneity of $Q$ to rewrite the 
left hand side of~(\ref{q.q^2.2}). 
Since $E'(3,4)=0$ for $q=1$, we may assume $q \ge 2$.
The Gordan series 
$\Gordan{Q}{Q}{Q}{q}{q}{q}{0}{0}{2}$ gives an identity
\[ (Q^2,Q)_2 = \frac{3 \, q-2}{2(2q-1)} \, Q \, (Q,Q)_2. \] 
(See \cite[Ch.~IV]{GrYo} for the derivation of the series.) 
We have shown that $Q$ divides $E'(3,4)$, hence $J$ is a polynomial 
covariant. 

We can continue the calculation of $J$ from 
(\ref{partialijf}), but it is easier to proceed as follows. By 
counting degrees, we see that $J(P,Q)$ is a joint covariant of 
$P,Q$ which is quadratic in $P,Q$ separately and has order $2p+2q-4$. 
We claim that every such joint covariant is a linear combination of 
\begin{equation} 
 Q^2 \, (P,P)_2, \quad P^2 \, (Q,Q)_2, \quad (P^2,Q^2)_2. 
\label{3cov} \end{equation}
This amounts to counting the number of copies of 
the representation $S_{2p+2q-4}$ inside 
$S_2(S_p) \otimes S_2(S_q)$. A straightforward expansion shows that 
there are three such copies (see \cite[\S 4.2]{Sturmfels}). 
It is easy to see by specialization that the covariants in 
(\ref{3cov}) are linearly independent for generic $P,Q$, so
they must form a basis for this space. 

Hence we may write $J$ as in (\ref{jpq}) for some constants $z_i$. 
Specialize to $P = x_0^p, Q = x_1^q$, then $(P,P)_2 = (Q,Q)_2 =0$ and 
$(P^2,Q^2)_2 = x_0^{2p-2} \, x_1^{2q-2}$. On the other hand, 
$Q^4 \, \He(U) = p \, q \, (p-q-1) \, x_0^{2p-2} \, x_1^{2q-2}$. 
This forces 
\[ z_3 = p \, q \, (p-q-1). \] 
Similarly, specialize $P,Q$ to the pairs 
$(x_0^p,x_0^{q-1} \, x_1)$ and $(x_0^{p-1} \, x_1,x_0^q)$, and 
get two more linear equations involving the $z_i$. Solving these, 
we get the theorem.  \qed 

This formula has a simple but interesting corollary. 
If $p=q+1$, then $\He \, (\frac{P}{Q})$ is identically zero. 

Finally write $\gC_e(F) = J(F^{2d-2e-2},He(F)^{d-e})$; then 
we can state a criterion which involves only polynomial 
covariants: 
\begin{Theorem} \sl Let $F$ be a binary $d$-ic. Then 
\[ F \in X^{(d-e,e)} \iff \gC_e(F)=\D(F)=0.\]
\end{Theorem} 

\subsection{A formula for $\D$} 
One can write down a formula for the covariant $\D$ in terms of 
compound transvectants. The proofs will only be sketched. Define 
\[ \begin{aligned} 
\xi_1 & = (2d-1)(2d-3)(2d-5)^3, \\ 
\xi_2 & = -9 \, (d-3)(2d-5)(2d-7)(2d-3)^2, \\ 
\xi_3 & = 4 \, (d-1)(d-3)(d-4)(2d-9)(4d-7). 
\end{aligned} \] 
\begin{Proposition} \sl 
If $d \ge 6$, then up to a multiplicative scalar
\begin{equation} \D(F) = \xi_1 \, (F^2,F^2)_6 + 
\xi_2 \, (F^2,(F,F)_2)_4 + \xi_3 \, (F^2,(F,F)_4)_2. 
\label{D(F).formula} \end{equation}
\end{Proposition} 
One may argue as follows: for $d \ge 6$, there are 
three copies of $S_{4d-12}$ inside $S_4(S_d)$, and a basis for this 
space is given by the three covariants which occur in 
(\ref{D(F).formula}). Hence $\D(F)$ can be written as their linear 
combination. To determine the actual coefficients, 
specialize to $F = x_0^{d-e} x_1^e, e=2,3$ (when $\D$ must vanish)
and solve a system of linear equations. \qed 

\begin{Remark} \rm
It is a priori clear that the $\xi_i$ should be rational functions 
in $d$ (or polynomials after clearing denominators). 
However, we can see no conceptual explanation of the fact that 
they should split into linear factors over ${\mathbf Q}$. 
\end{Remark}

These are the formulae in low degrees: 
\[ \D(F) = \begin{cases}  
(F^2,F^2)_6 & \text{for $d=3$,} \\ 
7 \, (F^2,F^2)_6 - 5 \, (F^2,(F,F)_2)_4 & \text{for $d=4$,} \\ 
129 \, (F^2,F^2)_6 - 250 \, (F^2,(F,F)_2)_4 & \text{for $d=5$.} 
\end{cases} \] 
To prove these, notice that there are two copies of 
$S_{4d-12}$ in $S_4(S_d)$ for $d =4,5$ and argue as before. For 
degree $3$ forms, $\D$ is simply the discriminant. 

\section{Binary forms-II} \label{binary.idealgen} 
In this section we write down the covariants which correspond to 
the quartic generators of $I_{X^{(d-e,e)}}$. 
Given a triple of integers $I=(i,j,k)$, define a covariant
\[ \E_I(F) = (((F,F)_{2i},F)_j,F)_k. \] 
Let $F$ be a general point of $X^{(d-e,e)}$, then we may write  
$F = x_0^{d-e} \, x_1^e$ after a change of variables. 
Using formula (\ref{expr.T}), 
\[ \E_I(F) = \omega_I \, x_0^{4(d-e)-(2i+j+k)} \, x_1^{4e-(2i+j+k)} \] 
where $\omega_I$ is the rational number 
\[ \begin{aligned} 
{} & \cN\left[\begin{array}{c} d-e,e\\ d-e,e\\ 2i \end{array} \right] 
\times \, 
\cN\left[\begin{array}{c} 2(d-e)-2i,2e-2i\\ d-e,e\\ j \end{array} \right] 
\times \\
& \cN\left[\begin{array}{c} 3(d-e)-(2i+j),3e-(2i+j) \\ d-e,e\\ k
\end{array} \right] \end{aligned} \]
Thus, as a monomial, $\E_I(F)$ depends only on the sum $2i+j+k$. Given 
triples  $I=(i,j,k),I'=(i',j',k')$ such that $2i+j+k = 2i'+j'+k'$, 
define 
\[ \Psi_{I,I'}(F) = \omega_I \, \E_{I'}(F) - \omega_{I'} \, \E_I(F). \] 
\begin{Proposition} \sl 
The locus $X^{(d-e,e)}$ is scheme-theoretically generated by 
the coefficients of all the covariants $\Psi_{I,I'}$. 
\end{Proposition} 
\demo The proof is in essense identical to 
\cite[Theorem 7.2]{AC1}, hence we omit the details. \qed 

\smallskip 

We have been unable to give `closed formulae' for the $\omega_I$, 
indeed this is directly traceable to the difficulty that no 
closed expression is known for a general Clebsch-Gordan coefficient
$C^{j_1,j_2,j}_{m_1,m_2,m}$.

\smallskip 

The covariants corresponding to the {\sl quadratic} generators 
of $I_{X^{(d-e,e)}}$ are easily described, they are 
$\{(F,F)_{2i}: e+1 \le i \le [\frac{d}{2}]\}$.

\section{Ternary quintics} \label{x32}
In this section we work out the case $n=2,(d-e,e)=(3,2)$, 
and describe the ideal generators invariant-theoretically. We have 
made rather heavy use of machine-computations, specifically the 
programs Macaulay-2 and Maple. 

\subsection{} Define generic forms 
\[ \begin{aligned}
L_1 & = a_0 \, x_0 + a_1 \, x_1 + a_2 \, x_2, \quad 
L_2 = b_0 \, x_0 + b_1 \, x_1 + b_2 \, x_2 \\ 
F & = c_0 \, x_0^5 + c_1 \, x_0^4 \, x_1 + \dots + c_{20} \, x_2^5, 
\end{aligned} \] 
where $a,b,c$ are independent indeterminates. 
Write $F = L_1^3 \, L_2^2$ and equate the 
coefficients of the monomials in $x_0,x_1,x_2$. This expresses 
each $c_i$ as a polynomial in $a_0,\dots,b_2$, and hence 
defines a ring map 
\[ \complex[c_0,\dots,c_{20}] \lra \complex[a_0,\dots,b_2]. \] 
The kernel of this map is $I_X$. We calculated it in Macaulay-2, and found 
that its resolution begins with 
\[ \dots 
\ra R(-4) \otimes M^{(4)} \oplus R(-3) \otimes M^{(3)} \ra 
R \ra R/I_X \ra 0, \] 
where $M^{(3)},M^{(4)}$ are vector spaces of dimensions $455$ and $1470$ 
respectively. Thus there are no generators in degrees $\ge 5$. 
Since $SL(V)$ acts on this resolution, 
the $M^{(i)}$ are $SL(V)$-modules.

\begin{Lemma} \sl We have the following isomorphisms of  
$SL(V)$-modules: 
\begin{equation} \begin{aligned} 
   M^{(3)} & = S_{(9,3)} \oplus S_{(9,0)} \oplus S_{(7,5)} \oplus S_{(7,2)} 
    \oplus S_{(6,3)} \oplus S_{(3,3)} \oplus S_{(3,0)}, \\ 
   M^{(4)} & = S_{(16,4)} \oplus S_{(14,6)} \oplus S_{(12,8)} 
           \oplus S_{(10,10)}. 
\end{aligned} \label{M3M4} \end{equation}
\end{Lemma} 

\demo 
Since $M^{(3)} = (I_X)_3$, the first isomorphism follows from Corollary 
\ref{gr}. (Throughout this example, all the inner and outer products of 
Schur functions were calculated using the Maple package `SF'.)

The degree $4$ piece of $I_X$ is a direct sum of two parts: 
multiples of degree $3$ generators by linear forms, and the new 
generators $M^{(4)}$. Hence 
\[ [(I_X)_4] = [M^{(3)} \otimes S_5] - [N^{(4)}] + [M^{(4)}], \] 
where $N^{(4)}$ denotes the module of first syzygies in 
degree $4$. (We know practically nothing about $N^{(4)}$, but we will 
see that this is no obstacle.) 
Now $[(I_X)_4]$ can be calculated by Corollary \ref{gr}, and 
$[M^{(3)} \otimes S_5]$ by the Littlewood-Richardson rule. Hence 
the difference $[M^{(4)}] - [N^{(4)}]$ is known, we write it as 
\[ [M^{(4)}] - [N^{(4)}] = {\mathcal Y} + {\mathcal Z}, \] 
where ${\mathcal Y}$ (resp.~${\mathcal Z}$) is a positive 
(resp.~negative) linear combination of Schur polynomials. 
The actual calculation shows that 
\[ {\mathcal Y} = [S_{(16,4)} \oplus S_{(14,6)} \oplus 
S_{(12,8)} \oplus S_{(10,10)}]. \] 
It follows that each summand on the right must appear in $M^{(4)}$. 
Now the direct sum has dimension 
$585 + 504 + 315 + 66 = 1470 = \dim M^{(4)}$. Hence $M^{(4)}$ must 
in fact coincide with this sum. \qed 

\subsection{} 
By the standard formalism of \cite{Littlewood}, 
a submodule $S_{(a,b)} \subseteq S_r(S_5)$ corresponds to a 
concomitant of degree $r$, order $a-b$ and class $b$ of ternary 
quintics. 

We will illustrate how to write down such a concomitant symbolically. 
For instance, let $\Psi$ correspond to the inclusion 
$S_{(16,4)} \subseteq M^{(4)}$. 
Decomposing $S_4(S_5)$, we detect that it has two copies of 
$S_{(16,4)}$, hence ternary quintics have two independent 
concomitants $\Psi_1,\Psi_2$ of degree $4$, order $12$ and class $4$. 
Now consider the following Young tableau of shape $(16,4)$ filled with 
four symbolic letters $\alpha,\beta,\gamma,\delta$, 
each occuring $5$ times: 
\[ \begin{array}{cccccccccccccccc} 
\alpha & \alpha & \alpha & \alpha & \alpha & \beta &  \beta &  \beta & 
\gamma & \gamma & \gamma & \gamma & \delta & \delta & \delta & \delta \\ 
\beta & \beta & \gamma & \delta 
\end{array} \] 
Reading this tableau {\sl columnwise}, we can construct the concomitant
\[ \Psi_1 = (\alpha \, \beta \, u)^2 \, (\alpha \, \gamma \, u) \, 
(\alpha \, \delta \, u) \; \alpha_x \; \beta_x^3 \; \gamma_x^4 \; 
\delta_x^4. \] 

We will abbreviate this as 
$\Psi_1= \langle 5,3,4,4| 0,2,1,1 \rangle$. 
(This means that in the top row of the tableau $\alpha$ occurs 
five times, followed by $\beta$ thrice etc. Such a notation is 
possible because we will choose all of our tableaux to be 
semistandard for the order $\alpha < \beta < \gamma < \delta$.)
Similarly, let $\Psi_2 = \langle 5,1,5,5| 0,4,0,0 \rangle$. 
To show that the $\Psi_i$ are linearly independent, it is sufficient to 
evaluate them on any specific form, in fact $F=x_0^5-x_1^5$ would do. 
(We checked this in Maple.) Hence $\Psi_1,\Psi_2$ form a basis 
of the space of concomitants of degree $4$, order $12$ and class $4$. 

Now write $\Psi = \eta_1 \, \Psi_1 + \eta_2 \, \Psi_2$, and evaluate on 
$x_0^3 \, x_1^2 \in X^{(3,2)}$. By hypothesis $\Psi$ must vanish 
identically, this gives the equation 
\[ (\frac{57}{2500} \, \eta_1 + \frac{3}{50} \, \eta_2) \, 
x_0^8 \, x_1^4 \, u_2^4=0. \] 
Hence $\eta_1:\eta_2 = 50:-19$, which determines $\Psi$ (up to 
a scalar). 

We have worked out the complete list for all the summands in (\ref{M3M4}). 
In degree $3$ (where we need only three symbolic letters), the concomitants 
are 
\begin{equation} \begin{array}{lll}
\langle5,4,1|0,1,3|0,0,1\rangle, & 
\langle5,3,3|0,2,0|0,0,2\rangle, & 
\langle5,3,0|0,2,4|0,0,1\rangle, \\
\langle5,3,1|0,2,2|0,0,2\rangle, & 
\langle5,2,1|0,3,2|0,0,2\rangle, & 
\langle5,1,0|0,4,2|0,0,3\rangle, \\ 
\langle5,1,1|0,4,0|0,0,4\rangle.
\end{array} \label{cov3} \end{equation} 
In degree $4$, they are 
\begin{equation} \begin{aligned}
50 \, & \langle 5,3,4,4|0,2,1,1 \rangle
- 19 \, \langle 5,1,5,5|0,4,0,0 \rangle, \\
5 \, & \langle 5,5,4,0|0,0,1,5 \rangle 
- 8 \, \langle 5,4,0,0|0,1,5,0 \rangle, \\
& \langle 5,1|0,4 \rangle^2 + 2 \, \langle 5,3,2,2| 0,2,3,3 \rangle,  \\ 
& \langle 5,3,2,0|0,2,3,5 \rangle.
\end{aligned} \label{cov4} \end{equation}  
In conclusion we have the following result: 
\begin{Proposition} \sl 
Let $F$ be a ternary quintic with zero scheme $C \subseteq \P^2$. Then 
$C$ consists of a triple line and a double line, iff all the concomitants 
in (\ref{cov3}) and (\ref{cov4}) vanish on $F$. 
\qed \end{Proposition} 

\bigskip 

\noindent {{\sc Acknowledgements:} \smaller 
The first author would like to express his gratitude to
Professors David Brydges and Joel Feldman for the invitation to 
visit the University of British Columbia.
The second author would like to thank Professor James Carrell for 
his invitation to visit UBC. We are indebted to Daniel Grayson and 
Michael Stillman (authors of Macaulay-2), John Stembridge 
(author of the `SF' package for Maple), 
the digital libraries maintained by the Universities of 
Cornell, G{\"o}ttingen and Michigan, as well as J-Stor and 
Project Gutenberg.} 

%%%%%%%%%%%%%%%%%%%%%%%%%%%%%%%%%%%%%%%%%%%%%%%%%%%%
\bibliographystyle{plain}

%%%%%%%%%%%%%%%%%%%%%%%%%%%%

\pagebreak 
\parbox{6cm}{\small 
{\sc Abdelmalek Abdesselam} \\ 
Department of Mathematics\\
University of British Columbia\\
1984 Mathematics Road \\
Vancouver, BC V6T 1Z2 \\ Canada. \\ 
{\tt abdessel@math.ubc.ca}} 
\hfill 
\parbox{5cm}{\small 
LAGA, Institut Galil\'ee \\ CNRS UMR 7539\\
Universit{\'e} Paris XIII\\
99 Avenue J.B. Cl{\'e}ment\\
F93430 Villetaneuse \\ France.} 

\vspace{1.5cm} 

\parbox{6cm}{\small 
{\sc Jaydeep Chipalkatti} \\ 
Department of Mathematics\\
University of Manitoba \\ 
433 Machray Hall \\ 
Winnipeg MB R3T 2N2 \\ Canada. \\ 
{\tt chipalka@cc.umanitoba.ca}}

\end{document}